\documentclass[preprint,sort&compress,12pt,draft]{elsarticle}
\usepackage{amsmath}
\usepackage{mathrsfs}
\usepackage{stmaryrd}
\usepackage{bm}
\usepackage{amsfonts}
\usepackage{amssymb}
\usepackage{amsthm}
\usepackage{xcolor}

\newcommand{\mm}{\mathrm}
\newcommand{\ml}{\mathcal}

\newcommand{\be}{\begin{equation}}
\newcommand{\bea}{\begin{equation}\begin{aligned}}
\newcommand{\beas}{\begin{equation*}\begin{aligned}}
\newcommand{\eeas}{\end{aligned}\end{equation*}}
\newcommand{\eea}{\end{aligned}\end{equation}}
\newcommand{\ee}{\end{equation}}

\begin{document}
\begin{frontmatter}
\title{Global Solutions and Asymptotic Behavior for the Three-dimensional\\ Viscous Non-resistive MHD System with Some Large Perturbations}

\author[fJ]{Youyi Zhao}
\ead{zhaoyouyi957@163.com}
\address[fJ]{School of Mathematics and Statistics, Fuzhou University, Fuzhou, 350108, China.}
\begin{abstract}
We revisit the global existence of solutions with some large perturbations to the incompressible, viscous, and non-resistive MHD system in a three-dimensional periodic domain, where the impressed magnetic field satisfies the Diophantine condition, and the intensity of the impressed magnetic field, denoted by $m$, is large compared to the  perturbations.
It was proved by Jiang--Jiang that the highest-order derivatives of the velocity increase with $m$,
and the convergence rate of the nonlinear system towards a linearized problem is of $m^{-1/2}$ in [F. Jiang and S. Jiang, Arch. Ration. Mech. Anal., 247 (2023), 96].
In this paper, we adopt a different approach by leveraging vorticity estimates to establish the highest-order energy estimate. This strategy prevents the appearance of terms that grow with $m$, and thus the increasing behavior of the highest-order derivatives of the velocity with respect to $m$ does not appear. Additionally, we use the vorticity estimate to demonstrate the convergence rate of the nonlinear system towards a linearized problem as time or $m$ approaches infinity. Notably, our analysis reveals that the convergence rate in $m$ is faster compared to the finding of  Jiang--Jiang. Finally, a key contribution of our work is the identification of an integrable time-decay of the lower dissipation, which can replace the time-decay of lower energy in closing the highest-order energy estimate. This finding significantly relaxes the regularity requirements for the initial perturbations.
\end{abstract}
\begin{keyword}
MHD fluids; incompressible fluids; global well-posedness; large initial perturbation; vorticity estimate; algebraic time-decay; convergence rate in term of the field intensity.
\end{keyword}
\end{frontmatter}


\newtheorem{thm}{Theorem}[section]
\newtheorem{lem}{Lemma}[section]
\newtheorem{pro}{Proposition}[section]
\newtheorem{cor}{Corollary}[section]
\newproof{pf}{\emph{Proof}}
\newdefinition{rem}{Remark}[section]
\newtheorem{definition}{Definition}[section]

\section{Introduction}\label{sec:01}

\numberwithin{equation}{section}
The dynamics of electrically conducting fluids interacting with magnetic field
can be described by the system of magnetohydrodynamics (MHD) equations.
In this paper, we shall investigate the following three-dimensional (3D) incompressible, viscous, and non-resistive MHD system in the domain $\Omega$:
\begin{equation}\label{0101}
\begin{cases}
\rho{v}_t+ \rho{v}\cdot\nabla {v}-\nu\Delta{v}+\nabla \left(p+{\lambda}|M|^2/2\right)
=\lambda M\cdot\nabla M,  \\[1mm]
{M}_t+{v}\cdot\nabla {M}={M}\cdot\nabla {v}, \\[1mm]
\mm{div}v=\mm{div}{M}=0.
\end{cases}
\end{equation}
We consider that the domain is a 3D periodic domain, i.e. $\Omega=\mathbb{T}^3$, where $\mathbb{T}:=\mathbb{R}/\mathbb{Z}$ the 1D-torus.
The unknowns ${v}:={v}({x},t)$, ${M}:={M}({x},t)$ and $p:=p({x},t)$
denote the velocity, the magnetic field and the kinetic pressure of MHD fluids, respectively;
the three positive (physical) parameters $\rho$, $\nu$ and $\lambda$
stand for the density, shear viscosity coefficient and permeability of vacuum dividing by $4\pi$, respectively.
\subsection{Background and motivation}
{The investigation of the MHD system has captured the attention of mathematicians over the past half a century. Many intriguing and challenging mathematical problems have been explored from both analytical and numerical simulation perspectives, leading to significant advancements. Below we review some relevant previous results that are connected to our research.

Based on the linearized non-resistive MHD system, Chandrasekhar first found that magnetic fields can inhibit the thermal instability in electrically conducting fluids \cite{CSHHSCPO}. Later numerical observations have further indicated that the energy in an inhomogeneous plasma dissipates at a rate independent of resistivity \cite{CFCCRI}. Consequently, it can be inferred that the (nonlinear) non-resistive MHD system may still exhibit dissipative behavior due to the magnetic field, thereby allowing for the existence of global solutions, at least for small initial data. This conclusion has been mathematically verified in existing literature, see \cite{BCSCSPLL, RXXWJHXZYZZF, TZWYJGw} for examples.
Moreover, in the context of the nonlinear MHD system, a sufficiently strong magnetic field is shown to reduce nonlinear interactions, inhibiting the formation of strong gradients \cite{RHK} and flow instabilities \cite{JFJSJMFMOSERT, JFJSSETEFP, JFJSOUI, WYTIVNMI}.

For the ideal MHD system, namely the system \eqref{0101} with $\mu=0$, Bardos--Sulem--Sulem employed the transformation of the Els\"asser variables and then adopted a hyperbolic approach to establish the global(-in-time) existence of classical solutions to the Cauchy problem of the 2D/3D ideal MHD system under the framework of H\"older space \cite{BCSCSPLL}. Their result showed that the system \eqref{0101} with $\mu=0$ is globally well-posed for the large initial perturbations when the intensity of the impressed magnetic field is sufficiently large (i.e. the condition of strong field intensity). It's worth noting that such global stability result is not expected for the 3D incompressible Euler equations.

A question naturally raises  whether the viscous and non-resistive MHD system \eqref{0101} also admits a global solution with large perturbation under strong field intensity. Later Zhang \cite{ZTGS} reformulated the magnetic field $M$ by using the stream function and then developed the vanishing mechanism of nonlinear terms based on the energy method and spectrum analysis to successfully establish the existence result of large perturbation solutions for the system \eqref{0101} under  the strong field intensity in a 2D whole space $\mathbb{R}^2$. Motivated by the magnetic inhibition theory in \cite{JFJSOMITIN}, Jiang--Jiang \cite{JFJSO2020} reformulated the analysis of the non-resistive MHD system in Lagrangian coordinates and also developed new vanishing mechanism of nonlinear terms to obtain similar large perturbation solutions for the system \eqref{0101} in a 2D periodic domain $\mathbb{T}^2$; moreover, they further established the convergence rate of the nonlinear system toward a linearized system as the intensity of  impressed magnetic fields goes to infinity. Recently, when the Diophantine condition was imposed on impressed magnetic fields, Jiang--Jiang further obtained a existence result and asymptotic behavior of global large perturbation solutions for the system \eqref{0101} in a 3D periodic domain $\mathbb{T}^3$ under Lagrangian coordinates \cite{JFJSO2021}. However, they need the extremely high regularity of initial data to close the process of energy estimates via the multi-layer energy method. Additionally, when attempting to control the highest-order derivative of velocity, a term involving the highest-order energy estimate increases with the field intensity $m$, see \cite[(3.33)]{JFJSO2021} for example.

In this paper, we substitutively exploit the vorticity estimate to establish the highest-order energy estimate. Such approach eliminates terms that would increase with   field intensity, avoiding the appearance of the increasing behavior of the highest-order derivatives of velocity with respect to field intensity $m$. This approach has the additional advantage of eliminating challenging terms hidden in the integral term involving the pressure. Additionally, since we don't need an extra step to balance field intensity-increasing terms (see \cite[(3.33) and (3.43)]{JFJSO2021}), we can reduce the regularity requirement of the initial data in \cite{JFJSO2021}.
Furthermore, we employ the vorticity estimate method to derive the convergence rate of the nonlinear system toward a linearized problem as time or field intensity goes to infinity. Our new convergence result indicates a faster convergence rate with respect to field intensity compared to the finding in \cite{JFJSO2021}. It's worth noting that the handling of the pressure term, as in \cite[(3.16)]{JFJSO2021}, also influences the convergence rate of the nonlinear system toward a linearized problem.

It is well-known that in Lagrangian coordinates, a non-resistive MHD fluid in a rest state can be viewed as composed of infinite fluid element lines that are parallel to the impressed magnetic field. These lines can be considered as elastic strings under magnetic tension, which have been discussed in \cite{JFJSOMITIN}. Furthermore, using Lagrangian coordinates allows us to mathematically circumvent many difficulties arising from a nonlinear term involving the magnetic field. Specifically, the magnetic tension term $\lambda M\cdot\nabla M$ becomes a linear term with dissipation along the direction of the given impressed magnetic field, expressed by  the directional derivative
 of flow function.
Since studying the non-resistive MHD system in Lagrangian coordinates is convenient and effective, as highlighted in \cite{JFJSCVPDE1, JFJSOUI, JFJSSETEFP, JFJSJMFMOSERT, WYTIVNMI, JFJSO2020, JFJSO2021}, we also follow the argument presented in \cite{JFJSO2021} to reformulate the MHD system \eqref{0101} into Lagrangian coordinates.}

\subsection{Reformulation in Lagrangian coordinates}\label{0608241825}
Let $(v,M)$ be the solution of the MHD system \eqref{0101}.
We use the flow map $\zeta$ to define the Jacobi matrix $\mathcal{A}:=(\mathcal{A}_{ij})_{3\times 3}$
via $\mathcal{A}^{\mm{T}}:=(\partial_j \zeta_i)^{-1}_{3\times 3}$, where the flow map $\zeta$ is the solution to the initial value problem
\begin{equation*}
\begin{cases}
 \partial_t\zeta (y,t)=v(\zeta(y,t),t) &\mbox{ in } \mathbb{T}^3\times(0,\infty),\\
\zeta (y,0)= \zeta^0(y) &\mbox{ in } \mathbb{T}^3.
\end{cases}
\end{equation*}
Here we have assumed that $\zeta^0:=\zeta^0(y):\mathbb{T}^3\to \mathbb{T}^3$ is diffeomorphism, and satisfies
the volume-preserving condition $\det\nabla \zeta^0=1$, and $"\det"$ denotes the determinant.
We denote the Eulerian coordinates by $(x,t)$ with $x=\zeta(y,t)$, whereas $(y,t)\in \mathbb{T}^3\times(0,\infty)$ stand for the
Lagrangian coordinates.
By the divergence-free condition $\mm{div}v=0$, the flow map $\zeta$ also satisfies the volume-preserving condition \cite{MAJBAL}:
\begin{equation}\label{202207011433n}
\det\nabla \zeta=1
\end{equation}
as well as $\det\nabla \zeta^0=1$.

Define the Lagrangian unknowns by
\begin{equation*}
\quad\quad(u, B, Q)(y,t)=(v, M, p+ \lambda|M|^2/2)(\zeta(y,t),t),\quad \mbox{for}\;(y,t)\in \mathbb{T}^3\times(0,\infty).
\end{equation*}
The equations for $(\zeta,u,q,B)$ in Lagrangian coordinates read as follows.
\begin{equation}\label{202109221247}
\begin{cases}
\zeta_t=u&\mbox{ in } \mathbb{T}^3 ,\\[1mm]
\rho u_t - \nu \Delta_{\mathcal{A}}u+\nabla_{\mathcal{A}}Q = \lambda  B\cdot \nabla_{\mathcal{A}} B
&\mbox{ in } \mathbb{T}^3 ,\\[1mm]
B_t-B\cdot\nabla_{\mathcal{A}}u=0
&\mbox{ in } \mathbb{T}^3 ,\\[1mm]
\mm{div}_{\mathcal{A}}u=\mm{div}_{\mathcal{A}}B=0 &\mbox{ in } \mathbb{T}^3
 ,\\[1mm]
( u, B)|_{t=0}=( u^0, B^0):=( v^0(\zeta^0), M^0(\zeta^0))
&\mbox{ in } \mathbb{T}^3.
\end{cases}
\end{equation}
Here we have written the differential operators $\nabla_{\mathcal{A}}$, $\mm{div}_\mathcal{A}$ with their actions given by
$(\nabla_{\mathcal{A}}f)_i:=\mathcal{A}_{ij}\partial_jf$, $\mm{div}_{\mathcal{A}}X:=\mathcal{A}_{ij}\partial_j X_{i}$
for appropriate scalar function $f$ and vector function $X$, and $\Delta_{\mathcal{A}}f:=\mm{div}_{\mathcal{A}}
\nabla_{\mathcal{A}}f$. It should be noted that the Einstein convention of summation over repeated indices has been used here, and $\partial_k:=\partial_{y_k}$.
Additionally, thanks to \eqref{202207011433n}, we have the following so called geometric identity:
\begin{equation}
\label{AklJ=0} \partial_l\mathcal{A}_{kl}=0.
 \end{equation}

We turn to study the equivalently MHD system \eqref{202109221247} in Lagrangian coordinates,
and will prove the global well-posedness of the problem \eqref{202109221247}
around the equilibrium state $(u, B)=(0, \bar{M})$ with
\begin{equation}\label{202107011602}
\bar{M}:=\varpi\omega,
\end{equation}
where $\omega\in\mathbb{R}^3$ and $\varpi\in\mathbb{R}^{+}$ denote the unit vector and the field intensity of $\bar{M}$, respectively;
and the unit vector $\omega\in\mathbb{R}^3$ satisfies the Diophantine condition \cite{CZZ2021,JFJSO2021}:
\begin{align}\label{202207041032}
\exists\;\mbox{a constant}\;c_{\omega}>0,\;\mbox{such that}\;|\chi\cdot\omega|\geqslant c_{\omega}|\chi|^{-3}
\;\;\mbox{for any}\;\chi\in\mathbb{Z}^3\setminus\{0\}.
\end{align}

Further, from  the differential version of magnetic flux conservation \cite{JFJSOMITIN}, we can see that
\begin{align}\label{01061800}
B=\bar{M}\cdot \nabla \zeta:=\partial_{\bar{M}}\zeta
\end{align}
provided the initial data $(\zeta^0,B^0)$ satisfies the frozen condition
\begin{align}&\label{abjlj0i}
B^0= \bar{M}\cdot \nabla \zeta^0:=\partial_{\bar{M}}\zeta^0.
\end{align}
\emph{Here and in what follows, the notation $f^0$ (or $f_0$) denotes the initial data of function $f(\cdot,t)$}.
It is worth noting that $B$ given by \eqref{01061800} automatically satisfies \eqref{202109221247}$_3$ and \eqref{202109221247}$_4$.
Finally, based on \eqref{01061800}, one can further calculate that
\begin{align}\label{0131}
B\cdot \nabla_{\mathcal{A}} B=(\bar{M}\cdot \nabla)^2\zeta=\varpi^2\partial_{\omega}^2\zeta.
\end{align}

Let
$$\mu=\nu/\rho,\;\; m^2= \lambda\varpi^2 /\rho,\;\;
q= Q /\rho,\;\; \eta=\zeta-y,\;\;\mbox{and}\;\; \eta^0=\zeta^0-y.$$
Consequently, under assumption \eqref{abjlj0i},
we  use  \eqref{0131} to transform \eqref{202109221247} into
\begin{equation}\label{202109221247nn}
\begin{cases}
\eta_t=u&\mbox{ in } \mathbb{T}^3 ,\\[1mm]
u_t +\nabla_{\mathcal{A}}q- \mu\Delta_{\mathcal{A}}u=m^2\partial_{\omega}^2\eta\quad
&\mbox{ in } \mathbb{T}^3 ,\\[1mm]
\mm{div}_{\mathcal{A}} u= 0
&\mbox{ in } \mathbb{T}^3 ,\\[1mm]
(\eta,u)|_{t=0}=(\eta^0, u^0)
&\mbox{ in } \mathbb{T}^3,
\end{cases}
\end{equation}
and
$$B=\varpi\partial_{\omega}(\eta+y)\;\;\mbox{and}\;\;\mathcal{A}=(\nabla \eta+I)^{\mm{T}},$$
where $I$ denotes the $3\times 3$ identity matrix.
The rest of this paper is devoted to demonstrating the global well-posedness for \eqref{202109221247nn}
with certain large perturbations corresponding to $m$.
For the sake of simplicity, we will continue to call $m$ the intensity of impressed magnetic field in the following discussions.

{The rest of this paper is organized as follows.
In Section \ref{main results}, we introduce the main results including the existence of a unique classical solution to the problem \eqref{202109221247nn} with a certain class of large initial perturbations corresponding to the   field intensity. Additionally, we discuss the convergence rate of the nonlinear system towards a linearized problem as time or $m$ tends to infinity, stated in Theorems \ref{thm:2022n} and \ref{thm:2022nn}.
In Sections \ref{sec:03N} and \ref{AB}, we provide detailed proofs of Theorems \ref{thm:2022n} and \ref{thm:2022nn}, respectively.
}

\section{Main result}\label{main results}

\subsection{Notations}

Before stating our result, we list some conventions for notation which will be frequently used in this paper.
Let
$\langle \cdot \rangle:=1+\cdot$,
$\int\cdot:=\int_{\mathbb{T}^3}\cdot$,
$(\cdot)_{\mathbb{T}^3}=\int\cdot \mathrm{d}y$.
We use $c>0$ to denote some ``universal" constants, which may depend on the domain but not on $m$, and may be different from line to line; $c_0>0$ denotes a "generic" constant independent of any parameters;
moreover
$a\lesssim b$ and $a\lesssim_0 b$ mean that $a\leqslant cb$ and $a\leqslant c_0b$, respectively.
Let $\alpha:=(\alpha_1, \alpha_2, \alpha_3)$ be  some multi-index of order $|\alpha|:=\alpha_1+\alpha_2+\alpha_3$.
$H^{k}:=H^{k}(\mathbb{T}^3)$ with non-negative integer $k$ being the usual Sobolev space, the norm of $H^{k}$ is denoted by $\|\cdot\|_{k}$,
$\underline{H}^{k}:=\{w\in H^{k}~|~(w)_{\mathbb{T}^3}=0\}$,
${H}_{1}^{k+1}:=\{w\in H^{k+1}~|~\mm{div}(\nabla w+I)=1\}$,
$\underline{H}_{1}^{k+1}:=\underline{H}^{k}\cap{H}_{1}^{k+1}$ and
$H^{k+2}_{*}:=\{w\in H^{k+2}~|~\zeta:=w+y:\mathbb{R}^3\rightarrow\mathbb{R}^3\;\mbox{is a $C^{k}$-diffeomorphism mapping}\}$.

In addition, we define a set of lower-order energies and dissipations: 
$$\begin{aligned}
&\mathcal{E}_{i}:=\|(\nabla\eta,u,m\partial_{\omega}\eta)\|_{i}^{2},\quad
\mathcal{D}_{i}:=\|(\nabla u,m\partial_{\omega}\eta)\|_{i}^{2}\quad\;\;\;\mbox{for}\;i=0,...,8.
\end{aligned}$$
Different from \cite{JFJSO2021},
the highst-order energy and dissipation are defined as follows.
$$\begin{aligned}
&\mathcal{E}_{H}:=\mathcal{E}_{12}=\|(\nabla\eta,u,m\partial_{\omega}\eta)\|_{12}^2,
\quad\mathcal{D}_{H}:=\mathcal{D}_{12}=\|(\nabla u,m\partial_{\omega}\eta)\|_{12}^2.
\end{aligned}$$

\subsection{ Global-in-time existence of large perturbation solutions}
The first main result of this paper is stated as follows.
\begin{thm}\label{thm:2022n}
Let $\bar{M}:=\varpi\omega$ and the unit vector $\omega\in\mathbb{R}^3$ satisfy the Diophantine condition \eqref{202207041032}.
There exist positive constants $c_1\geqslant4$, $c_2>0$ and a sufficiently small constant $c_3\in(0,1]$ such that,
for any initial data $(\eta^0, u^0)\in (H^{13}_{*}\cap\underline{H}_{1}^{13})\times \underline{H}^{12}$ and $m$
satisfying $\mm{div}_{\mathcal{A}^0}u^0=0$
and the condition of strong magnetic field
\begin{align}\label{202207041654}
\frac{\,\;\max\{(c_1\mathcal{E}_{H}^0e^{c_2\Psi})^{1/2}, c_1\mathcal{E}_{H}^0e^{c_2\Psi}\}\;}{m}\leqslant c_3,
\end{align}
the MHD problem \eqref{202109221247nn} admits a unique solution
$(\eta, u, q)\in (H^{13}_{*}\cap\underline{H}_{1}^{13})\times \underline{H}^{12}\times\underline{H}^{12}$ on $[0,\infty)$. Moreover, the solution $(\eta,u)$ enjoys that for $t\geqslant0$,
\begin{enumerate}
\item[(1)] the stability estimate:
\begin{align}
&\label{202207041516n}
\mathcal{E}_{i}(t)+\int_0^{t}\mathcal{D}_{i}(\tau)\mm{d}\tau\lesssim\mathcal{E}_{i}^0\quad\mbox{ for }0\leqslant i\leqslant8,\\[1mm]
&\label{202207041516}
\mathcal{E}_{H}(t)\lesssim\mathcal{E}_{H}^0 e^{c_2\Psi},\\[1mm]
&\label{202207041516nn}
\mathcal{E}_{H}(t)+\int_0^{t}\mathcal{D}_{H}(\tau)\mm{d}\tau\lesssim\mathcal{E}_{H}^0(1+\Psi e^{c_2\Psi}):=\ml{E}_{\mm{total}}^{0},
\end{align}
where $\Psi$ is defined in \eqref{202210141600};
\item[(2)] the decay-in-time estimate:
\begin{align}
&\label{202109090810nnm}
\langle t\rangle^{5/4}\mathcal{E}_{3}(t)+\int_0^{t}\langle \tau\rangle^{\frac{(1+5/4)}{2}}\mathcal{D}_{3}(\tau)\mm{d}\tau\lesssim\Xi,\\
&\label{202109090810}
\sum_{j=0}^{2}\langle m^{-1}\rangle^{2j}\left(\langle t\rangle^{(2-j)}\mathcal{E}_{4j}(t)+\int_0^{t}\langle \tau\rangle^{(2-j)}\mathcal{D}_{4j}(\tau)\mm{d}\tau\right)
\lesssim \Xi,
\end{align}
where
\begin{align}\label{202210141602}
\Xi:=\sum_{j=0}^{2}\langle m^{-1}\rangle^{2j}\mathcal{E}_{4j}^0.
\end{align}
\end{enumerate}
In addition, the solution $(\eta,u, q)$ satisfies the additional estimates \eqref{al1}--\eqref{alvarpi1}, \eqref{qestimate} and
\begin{align}\label{202210061225}
\sup_{t\in[0,\infty)}\|\eta(t)\|_{11}\lesssim 1.
\end{align}
\end{thm}

\begin{rem}
In Theorem \ref{thm:2022n}, it have been assumed that $(\eta^0)_{\mathbb{T}^3}=(u^0)_{\mathbb{T}^3}=0$.
Otherwise, we can define $\bar{\eta}^0=\eta^0-(\eta^0)_{\mathbb{T}^3}$ and $\bar{u}^0=u^0-(u^0)_{\mathbb{T}^3}$
such that $(\bar{\eta}^0)_{\mathbb{T}^3}=(\bar{u}^0)_{\mathbb{T}^3}=0$.
Consequently, by virtue of Theorem \ref{thm:2022n}, there exists a unique classical solution $(\bar{\eta},\bar{u}, \bar{q})$
to the problem \eqref{202109221247nn} with initial data $(\bar{\eta}^0,\bar{u}^0)$ in place of $(\eta^0,u^0)$.
Therefore, it is easy to verify that $(\eta,u,q):=(\bar{\eta}+t(u)_{\mathbb{T}^3}+(\eta^0)_{\mathbb{T}^3},\bar{u}+(u)_{\mathbb{T}^3}, \bar{q})$ is the unique classical solution
of the problem \eqref{202109221247nn} with initial data $(\eta^0,u^0)$.
\end{rem}
\begin{rem}
Recalling the derivation of the decay-in-time estimate \eqref{202109090810},
it is easy to see that the higher the regularity of the initial data is, the faster the decay in time of low-order derivatives of solutions is. Consequently, we can establish a result for a faster decay-in-time estimate, similar to the finding in \cite{JFJSO2021}. In fact, we can  further  establish a result for an almost exponential decay-in-time estimate, as demonstrated in \cite{TZWYJGw}.
\end{rem}

Now we briefly sketch the proof of Theorem \ref{thm:2022n}.
Let us first recall the basic energy identity of the initial value problem \eqref{202109221247nn} as follows.
$$\frac{1}{2}\frac{\mm{d}}{\mm{d}t}\left(\|u\|^2_0+ \|m\partial_{\omega} \eta\|_{0}^2 \right)+ \mu\|\nabla_{\ml{A}} u\|_{0}^2=0.$$
Integrating the above identity over $(0,t)$ yields that
\begin{align}\label{202210241617}
\|u(t)\|^2_0+ \|m\partial_{\omega} \eta(t)\|_{0}^2+2\mu\int_{0}^{t}\|\nabla_{\ml{A}} u(\tau)\|_{0}^2\mm{d}\tau
=\|u^0\|^2_0+ \|m\partial_{\omega} \eta^0\|_{0}^2:=I^0.
\end{align}
In particular, we have
\begin{align}\label{202210241621}
\|\partial_{\omega} \eta(t)\|_{0}^2\leqslant I^0/m^2.
\end{align}
We call $I^0$ the initial data mechanical energy.
It is easy to see that, for given (or even sufficiently large) initial data mechanical energy $I^0$,
$$\|\partial_{\omega} \eta(t)\|_{0}\rightarrow 0\quad\mbox{as}\;\;m\rightarrow\infty.$$
The above fundamental relation leads us to anticipate that the deformation quantity $\nabla\eta$ may be small when $m$ is sufficiently large. Fortunately, in the relatively high-regularity framework, this expectation holds true when the direction of the impressed magnetic field $\omega$ satisfies the Diophantine condition \eqref{202207041032}. This condition ensures a generalized Poincar\'e's inequality \eqref{poincareg}, providing the necessary  dissipation estimates of $\eta$ from $\partial_{\omega}^2\eta$.

This key observation motivates us to formally conclude that the problem \eqref{202109221247nn} may be approximated by the corresponding linearized problem for sufficiently large $m$. Since it is easy to prove that the linear problem has global solutions with a large initial perturbation, we could expect that the initial value problem \eqref{202109221247nn} may also admit a global large solution corresponding to the field intensity.
Our basic ingredient of establishing the existence of \eqref{202109221247nn} is the natural energy method, which comes from testing the momentum equations by $u$ and $\eta$:
\begin{align*}
&\frac{1}{2}\frac{\mm{d}}{\mm{d}t}\left(\|u\|_0^2+\|m\partial_{\omega}\eta\|_0^2\right)
+\|\sqrt{\mu}\nabla u\|_0^2\lesssim \mathcal{N}^{1},\\[1mm]
&\frac{1}{2}\frac{\mm{d}}{\mm{d}t}\|\sqrt{\mu}\nabla\eta\|_{0}^2+\int\eta\cdot u_{t}\mm{d}y
+\|m\partial_{\omega}\eta\|_0^2\lesssim \mathcal{N}^{2},
\end{align*}
where $\mathcal{N}^{1}$ and $\mathcal{N}^{2}$ represent integrals involving the nonlinear terms.
This demonstrates that
the dissipation still provides no direct control of the energy estimate of $\nabla\eta$ in $L^2(\Omega)$ due to strong coupling between the velocity and the magnetic field, as well as the weak dissipation from the magnetic tension,
and thereby leads to the presence of some nonlinear terms that cannot be controlled by dissipation alone.
The same type of problem is also encountered in the high-regularity energy estimates.

It should be noted that we need to work within a relatively high-regularity framework due to the usage of the generalized Poincar\'e's inequality \eqref{poincareg}. But difficulties thereby arise from nonlinear terms due to the 3D structure of the motion equations and the loss of derivatives resulting from the usage of the generalized Poincar\'e's inequality \eqref{poincareg}.
It is worth noting that similar difficulties also arise in the study of surface wave problems of viscous pure fluids
where the dissipation is weaker than the energy as well.
To overcome this difficulty, Guo--Tice introduced the celebrated two-layer energy method to consider
separately the low-order energy  with decay-in-time and the bounded high-order energy  \cite{GYTIDAP,GYTIAE2}.
In the spirit of Guo--Tice's method,
below we introduce a new version of the three-tier energy method. This approach considers separately the low-order energy/dissipation with time-decay/integrable time-decay and the bounded high-order energy to overcome the difficulties mentioned above. Additionally, we perform a vorticity energy estimate to establish the highest-order energy, preventing the increasing behavior of the highest-order energy concerning the field intensity $m$.
As a result, we can prove the existence of large perturbation global solutions to the problem \eqref{202109221247nn}, as stated in Theorem \ref{thm:2022n}.

Note that the key step to prove Theorem \ref{thm:2022n} is to derive the a \emph{priori} estimates \eqref{202207041516n}--\eqref{202207041516nn} with respect to sufficiently large $m$. By the above analysis, we shall perform three-tier energy method to prove that
there exist two constants $K$ and $\delta$, such that
\begin{align}\label{pr}
\sup_{t\in[0,T]}\left(\|\eta\|_{13}^2+\|(u,m\partial_{\omega}\eta)\|_{12}^2\right)\leqslant K^2/4,
\end{align}
provided that
\begin{align}\label{prio1}
\sup_{t\in[0,T]}\left(\|\eta\|_{13}^2+\|(u,m\partial_{\omega}\eta)\|_{12}^2\right)\leqslant K^2\quad\mbox{for any given}\;\; T>0
\end{align}
and
\begin{align}\label{prio2}
\max\{K, K^2\}/m\in(0,\delta],\;\;\delta\;\mbox{is sufficiently small}.
\end{align}
Roughly speaking, we designate \eqref{202207041516n} with $0\leqslant i\leqslant7$, \eqref{202207041516n} with $i=8$, and \eqref{202207041516} as the lower-order, higher-order, and highest-order energy inequalities, respectively. Under the assumptions of \eqref{prio1} and \eqref{prio2}, we first establish the (uniform in time and in the field intensity $m$) lower-order and higher-order energy inequalities, and then we obtain the lower-order energy with decay-in-time in \eqref{202109090810}. With \eqref{202109090810} established in hand, we further deduce an integrable time-decay of the lower dissipation \eqref{202109090810nnm}. Leveraging the  lower-order energy with decay-in-time in \eqref{202109090810} and the integrable time-decay of the lower dissipation \eqref{202109090810nnm}, we eventually close the highest-order energy inequality \eqref{202207041516}, yielding \eqref{pr}, please refer to Subsection \ref{stabilitye} for a detailed derivation.

It's worth noting that the constant $K$ will be determined later by \eqref{202209171746}, the smallness of $\delta$ may depend solely on the parameter $\mu$ and the unit vector $\omega$, and we have performed a vorticity energy estimate to establish the highest-order energy. Based on these considerations, and coupled with the existence of a unique local solution and a continuity argument as in \cite{JFJSO2021}, we immediately obtain Theorem \ref{thm:2022n}. The detailed proof will be presented in Section \ref{sec:03N}.

Finally, we remark the main differences between the proof of Theorem \ref{thm:2022n} and the stability result in \cite{JFJSO2021}:
\begin{enumerate}
\item[(1)]
{
In \cite{JFJSO2021}, the authors directly applied $\partial^{\alpha}$ to \eqref{202109221247nn}$_2$
to derive the highest-order energy estimates, resulting in the increasing behavior of the highest-order derivatives of velocity with respect to the magnetic field intensity $m$. This occurs because the $m$-increasing terms are hidden in the integral term $\int\partial^{\alpha}\nabla q\cdot\partial^{\alpha}u\mm{d}y=
-\int\partial^{\alpha} q\partial^{\alpha}\mm{div}u\mm{d}y$.
To address this issue, we instead employ the vorticity estimate method.
By using the vorticity estimate method, we can establish highest-order energy estimates without the undesirable $m$-increasing behavior. Consequently, both the low-order and highest-order energy evolutions in the proof process become much simpler compared to those in \cite{JFJSO2021}.
}
\item[(2)]
{In \cite{JFJSO2021},
the authors utilized the three-tier energy method, which necessitates the time-decay of lower-order energy $\mathcal{E}_{3}$ at a rate of $\langle t\rangle^{-(1+s)}$ with $s>1$.
Consequently, they had to establish the following lower-order energy with decay-in-time
\begin{align*}
\sum_{i=0}^{3}\langle m^{-1}\rangle^{2i}\left(\langle t\rangle^{(3-i)}\mathcal{E}_{4i}(t)+\int_0^{t}\langle \tau\rangle^{(3-i)}\mathcal{D}_{4i}(\tau)\mm{d}\tau\right)
\lesssim\sum_{i=0}^{3}\langle m^{-1}\rangle^{2i}\mathcal{E}_{4i}^0
\end{align*}
to achieve the required time-decay rate for $\mathcal{E}_{3}$.
In this paper, we find that an integrable time-decay of the lower dissipation $\mathcal{D}_{3}$ at a rate of $\langle t\rangle^{-(1+s)}$ with $s>0$ (see \eqref{202109090810nnm}--\eqref{202109090810}) can be employed instead of the time-decay of lower energy $\mathcal{E}_{3}$ at a rate of $\langle t\rangle^{-(1+s)}$ with $s>1$ to close the highest-order energy estimate.
This substitution significantly reduces the regularity requirements of the initial data concerning the fourth-order spatial derivative, please refer to \eqref{202207041516}--\eqref{202207041516nn} for details.
}
\end{enumerate}

\subsection{Vanishing phenomena of the nonlinear interactions}
First of all, we state the global existence of large perturbation solutions for the corresponding linear problem of \eqref{202109221247nn}.
\begin{pro}\label{2022nnnn}
Let the initial date $(\eta^0,u^0)$ be given in Theorem \ref{thm:2022n}.
There exists a function pair $(\eta^{\mm{r}},u^{\mm{r}})\in \underline{H}^{13}\times\underline{H}^{12}$ such that,
the following linear pressureless problem
\begin{equation}\label{202109221115}
\begin{cases}
\eta^{L}_t=u^{L}&\mbox{ in } \mathbb{T}^3 ,\\[1mm]
u^{L}_t - \mu\Delta u^{L}=m^2\partial_{\omega}^2\eta^{L}
&\mbox{ in } \mathbb{T}^3 ,\\[1mm]
\mm{div} u^{L}= 0
&\mbox{ in } \mathbb{T}^3 ,\\[1mm]
(\eta^{L}, u^{L})|_{t=0}=(\eta^0+\eta^{\mm{r}}, u^0+u^{\mm{r}})
&\mbox{ in } \mathbb{T}^3,
\end{cases}
\end{equation}
admits a unique global solution $(\eta^{L},u^{L})$ on $[0,\infty)$.
Moreover, the solution $(\eta^{L},u^{L})$ enjoys that for $t\geqslant0$,
$$(\eta^{L})_{\mathbb{T}^3}=(u^{L})_{\mathbb{T}^3}=\mm{div}\eta^{L}=\mm{div}u^{L}=0$$
and
\begin{align}
&\label{20221014112n}
\mathcal{E}_{12}^{L}(t)+\int_0^{t}\mathcal{D}_{12}^{L}(\tau)\mm{d}\tau\lesssim\mathcal{E}_{12}^0(1+\mathcal{E}_{3}^0)
\lesssim\ml{E}_{\mm{total}}^0,\\[1mm]
&\label{20221014112nN}
\mathcal{E}_{j}^{L}(t)+\int_0^{t}\mathcal{D}_{j}^{L}(\tau)\mm{d}\tau\lesssim\mathcal{E}_{j}^0
\;\;\quad{for}\;\;0\leqslant j\leqslant11,\\[1mm]
&\label{202210141123}
\sum_{j=0}^{2}\langle m^{-1}\rangle^{2j}\left(\langle t\rangle^{(2-j)}\mathcal{E}_{4j}^{L}(t)+\int_0^{t}\langle \tau\rangle^{(2-j)}\mathcal{D}_{4j}^{L}(\tau)\mm{d}\tau\right)
\lesssim\Xi,
\end{align}
where
$$\mathcal{E}_{j}^{L}:=\|(\nabla\eta^{L}, u^{L}, m\partial_{\omega}\eta^{L})\|_{j}^2\quad\mbox{and}\quad
\mathcal{D}_{j}^{L}:=\|(\nabla u^{L}, m\partial_{\omega}\eta^{L})\|_{j}^2.$$

In addition, the function pair $(\eta^{\mm{r}},u^{\mm{r}})$ satisfies
\begin{align}
&\label{202210301600}
\mm{div}(\eta^0+\eta^{\mm{r}})=\mm{div}( u^0+u^{\mm{r}})=0,\\[1mm]
&\label{202210301601}
\|\eta^{\mm{r}}\|_{l}\lesssim\|\mm{div}\eta^{0}\|_{l-1}\lesssim\|\eta^0\|_{3}\|\eta^0\|_{l},\\[1mm]
&\label{202210301602}
\|u^{\mm{r}}\|_{k}\lesssim\|\mm{div}_{\tilde{\ml{A}}^0}u^0\|_{k-1}\lesssim\|\eta^0\|_{3}\|u^0\|_{k}+\|\eta^0\|_{k}\|u^0\|_{3},\\[1mm]
&\label{202210301600n}
\|\partial_{\omega}\eta^{\mm{r}}\|_{k}\lesssim\|\partial_{\omega}\mm{div}\eta^0\|_{k-1}
\lesssim\|\eta^0\|_{3}\|\partial_{\omega}\eta^0\|_{k}+\|\eta^0\|_{k}\|\partial_{\omega}\eta^0\|_{3},
\end{align}
where $1\leqslant k\leqslant 12$ and $1\leqslant l\leqslant 13$.
\end{pro}
\begin{pf}
We can refer to \cite[Section 4]{JFJSO2021} for the  proof of Proposition \ref{2022nnnn}.
\hfill$\Box$
\end{pf}

With Theorem \ref{thm:2022n} and Proposition \ref{2022nnnn} in hand,
now we state the vanishing phenomena of the nonlinear interactions with respect to the time $t$ and the field intensity $m$.
\begin{thm}\label{thm:2022nn}
Let $(\eta,u,q)$ be the global solution of the initial value problem \eqref{202109221247nn} given in Theorem \ref{thm:2022n},
$(\eta^{L},u^{L})$ be the global solution of the linear problem \eqref{202109221115} constructed in Proposition \ref{2022nnnn},
denote $(\eta^{d}, u^{d}):=(\eta,u)-(\eta^{L},u^{L})$, then for any $t\geqslant0$,
\begin{itemize}
\item[(1)] asymptotic behavior with respect to $m$:
\begin{align}
&\label{202210011616nnnn}
\mathfrak{E}_{i+1}^{d}(t)+\int_0^{t}\mathfrak{D}_{i+1}^{d}(\tau)\mm{d}\tau\lesssim m^{-2}\Gamma
,\;\;\mbox{for}\;\;0\leqslant i\leqslant7,\\[1mm]
&\label{202210011616}
\mathfrak{E}_{9}^{d}(t)+\int_0^{t}\mathfrak{D}_{9}^{d}(\tau)\mm{d}\tau
\lesssim m^{-\frac{3}{2}}
\Pi;
\end{align}
\item[(2)] asymptotic behavior with respect to $m$ and time $t$:
\begin{align}
&\label{202210011616nn}
\sum_{j=0}^{2}\langle m^{-1}\rangle^{2j}\left(\langle t\rangle^{(2-j)}\mathfrak{E}_{4j+1}^{d}(t)
+\int_0^{t}\langle \tau\rangle^{(2-j)}\mathfrak{D}^{d}_{4j+1}(\tau)\mm{d}\tau\right)
\lesssim m^{-\frac{3}{2}}\Upsilon,
\end{align}
where
\begin{align}
&\mathfrak{E}_{i}^{d}:=\|(\nabla\eta^{d}, u^{d}, m\partial_{\omega}\eta^{d})\|_{i}^2\quad\mbox{and}\quad
\mathfrak{D}_{i}^{d}:=\|(\nabla u^{d}, m\partial_{\omega}\eta^{d})\|_{i}^2,\quad\;1\leqslant i\leqslant9,\nonumber\\[2mm]
&\Gamma:=\big(\sqrt{\ml{E}_{\mm{total}}^0}+\ml{E}_{\mm{total}}^0\big)^3,\nonumber\\[1mm]
&\Pi:=\max\{1, m^{-1/4},m^{-1/2}\}\big(\sqrt{\ml{E}_{\mm{total}}^0}+\ml{E}_{\mm{total}}^0\big)^3,\nonumber\\[1mm]
&\Upsilon:=\langle m^{-1}\rangle^4\max\{1, m^{-1/4},m^{-1/2}\}\big(\sqrt{\ml{E}_{\mm{total}}^0}+\ml{E}_{\mm{total}}^0\big)^2
\big(\sqrt{\ml{E}_{\mm{total}}^0}+\ml{E}_{\mm{total}}^0+\Xi\big).\nonumber
\end{align}
\end{itemize}
\end{thm}
\begin{rem}
Here we remark the main differences in the proof of Theorem \ref{thm:2022nn} and the asymptotic results in \cite{JFJSO2021}.
In \cite{JFJSO2021}, the authors directly applied $\partial^{\alpha}$ to \eqref{20210922111nn}$_2$ to derive the error energy estimates for $(\eta^{d}, u^{d})$, leading to a subtle term $\int\partial^{\alpha}\nabla q\cdot\partial^{\alpha}u^{d}\mm{d}y=-\int\partial^{\alpha}q\partial^{\alpha}\mm{div}u^{d}\mm{d}y$ as similar as the highest-order energy evolution in \cite{JFJSO2021}.
To achieve the convergence rate with respect to $m$, they had to use \eqref{20210922111nn}$_3$ to replace $\mm{div}u^{d}$, and consequently, whether in the higher-order error energy or in the lower-order error energy, the convergence rate of the error energy is  just $m^{-1}$.
To improve the the convergence rate of the nonlinear system toward a linearized problem, we instead employ the $\mm{curl}$-estimate method.
Our asymptotic result  demonstrates that the convergence rate of
the error energy $\mathfrak{E}_{i+1}^{d}$ with $0\leqslant i\leqslant7$ is   $m^{-2}$,
and thus the convergence rate of the error energy $\mathfrak{E}_{9}^{d}$ is  $m^{-3/2}$.
Therefore, the convergence rate with respect to $m$ in Theorem \ref{thm:2022nn} is faster than in \cite{JFJSO2021}, please refer to Section \ref{AB} for details.
\end{rem}
\begin{rem}  In addition, for $m\geqslant1$, the estimates \eqref{202210011616}--\eqref{202210011616nn} can be simplified as follows.
\begin{align}
&\label{2022103n1610}
\mathfrak{E}_{9}^{d}(t)+\int_0^{t}\mathfrak{D}_{9}^{d}(\tau)\mm{d}\tau
\lesssim m^{-3/2}\Gamma,\\[1mm]
&\label{2022103n1606}
\sum_{j=0}^{2}\left(\langle t\rangle^{(2-j)}\mathfrak{E}_{4j+1}^{d}(t)
+\int_0^{t}\langle \tau\rangle^{(2-j)}\mathfrak{D}^{d}_{4j+1}(\tau)\mm{d}\tau\right)
\lesssim m^{-3/2}\Gamma.
\end{align}
\end{rem}

\section{Proof of Theorem \ref{thm:2022n}}\label{sec:03N}
This section is dedicated to the proof of Theorem \ref{thm:2022n}.
We begin by listing some well-known mathematical results and introducing some preliminary estimates, which will be utilized in the derivation of energy evolution. Subsequently, we derive the energy estimates of  solutions.
Let $(\eta,u,q)$ be a solution of the  problem \eqref{202109221247nn} with the initial data $(\eta^0,u^0)\in\underline{H}_{1}^{13}\times\underline{H}^{12}$ and $\mm{div}_{\ml{A}^0}u^0=0$.
It is evident that $(\eta, u)$ automatically satisfies $(\eta)_{\mathbb{T}^3}=(u)_{\mathbb{T}^3}=0$.
Furthermore, we assume that the solution $(\eta,u,q)$ and $K$ satisfy $(q)_{\mathbb{T}^3}=0$ and
\eqref{prio1}--\eqref{prio2}, where $K$ will be determined later by \eqref{202209171746}.

\subsection{Well-known estimates and inequalities}
Now we list some well-known mathematical results, which will be used in the rest sections.
\begin{lem}\label{201806171834}
We have the following basic inequalities:
\begin{enumerate}[(1)]
\item
Poincar\'e's inequality (see \cite[Lemma 1.43]{NASII04}):
\begin{align}
\label{poincare}
\|w\|_{0}\lesssim \|\nabla w\|_{0}\;\quad\mbox{for}\;w\in\underline{H}^1.
\end{align}
\item
Generalized Poincar\'e's inequality (see \cite{CZZ2021}):
if the unit vector $\omega\in\mathbb{R}^3$ satisfies the Diophantine condition \eqref{202207041032}, it holds that
\begin{align}
\label{poincareg}
\|w\|_{0}\lesssim \|\partial_{\omega} w\|_{3}\;\quad\mbox{for}\;w\in\underline{H}^4.
\end{align}
 \item Interpolation inequalities (see \cite{ARAJJFF}): for any given $0\leqslant j<i$,
\begin{align}
&\label{interpolation2}
\|w\|_{j}\lesssim\|w\|_{0}^{1-j/i}\|w\|_{i}^{j/i}\;\quad\mbox{for}\;w\in{H}^{i}.
\end{align}
 \item Product estimates (see Section 4.1 in \cite{JFJSCVPDE1}): for any given $i\geqslant0$,
\begin{equation}
\label{product21}
\|fg\|_{i}\lesssim
\begin{cases}
\|f\|_{2}\|g\|_{i}
&\;\;\mbox{for}\;f,g\in{H}^{i}\cap H^{2},\;\;i=0,1;\\[1mm]
\|f\|_{2}\|g\|_{i}+\|f\|_{i}\|g\|_{2}
&\;\;\mbox{for}\;f,g\in{H}^{i}\cap H^{2},\;\;i\geqslant2.
\end{cases}
\end{equation}
                    \end{enumerate}
\end{lem}

\begin{lem}\label{lem:1642}
Hodge-type elliptic estimate: Let $w\in H^{i}$ with $i\geqslant1$, then it holds that
\begin{align*}
\|w\|_{i}\lesssim\|w\|_{0}+\|\mm{curl}w\|_{i-1} +\|\mm{div}w\|_{i-1}.
\end{align*}
In particular, if $w$ further satisfies zero-average condition $(w)_{\Omega}=0$, then
\begin{align}\label{hodgeellipticnn}
\|w\|_{i}\lesssim\|\mm{curl}w\|_{i-1} +\|\mm{div}w\|_{i-1}.
\end{align}
\end{lem}
\begin{pf}
The above two estimates are well-known and follow from the following identity: $$-\Delta w=\mm{curl}\mm{curl}w-\nabla\mm{div}w.$$
\hfill$\Box$
\end{pf}

\subsection{Preliminary estimates}
This subsection is devoted to the introduction of   preliminary estimates.
\begin{lem}\label{lem220914}
Let $\eta$ satisfy
\begin{align}\label{small1}
\sup_{t\in[0,T]}\|\eta\|_3\lesssim_{0}\delta\in(0,1],
\end{align}
then there holds that
\begin{align}
&\label{al1}
\|\tilde{\ml{A}}\|_{k}\lesssim_{0}\|\eta\|_{k+1},\;\;\hbox{ for }\;0\leqslant k\leqslant 12,\\[1mm]
&\label{alt1}
\|{\ml{A}}_{t}\|_{k}\lesssim_{0}
\begin{cases}
\|u\|_{i+1}
& \;\hbox{ for }\; 0\leqslant k\leqslant2;\\[1mm]
\|u\|_{i+1}+\|\eta\|_{i+1}\|u\|_{3}
& \;\hbox{ for }\; 3\leqslant k\leqslant12,
\end{cases}\\
&\label{alvarpi1}
\|\partial_{\omega}{\ml{A}}\|_{k}\lesssim_{0}
\begin{cases}
\|\partial_{\omega}\eta\|_{i+1}
& \;\hbox{ for }\; 0\leqslant k\leqslant2;\\[1mm]
\|\partial_{\omega}\eta\|_{i+1}+\|\eta\|_{i+1}\|\partial_{\omega}\eta\|_{3}
& \;\hbox{ for }\; 3\leqslant k\leqslant11.
\end{cases}
\end{align}
Here and in what follows $\tilde{\ml{A}}:={\ml{A}}-I$.
\end{lem}
\begin{pf}
One can refer to \cite[Lemma 3.1]{JFJSO2021} for the proof of Lemma \ref{lem220914}.
Here, we provide a brief sketch.
Recalling that \eqref{202207011433n} and $\zeta:=\eta+y$, then
\begin{align*}
\mm{det}\left(\nabla\eta+I\right)=1.
\end{align*}
It follows from the definition of matrix $\ml{A}$ that
the explicit form of the entries of $\tilde{\mathcal{A}}:=\ml{A}-I$ can be expressed as the quadratic polynomial $P_{2}(\partial_{i}\eta_{j})$
with respect to $\partial_{i}\eta_{j}$, where $P_{2}(\cdot)$ satisfies $P_2(0)=0$.
Moreover, it is easy to see that
$${\ml{A}}_{t}=\tilde{\ml{A}}_{t},\;\;\mbox{and}\;\;\partial_{\omega}\ml{A}=\partial_{\omega}\tilde{\ml{A}}.$$
Consequently, by applying $\|\cdot\|_{k}$, $\|\partial_{t}\cdot\|_{k}$, and $\|\partial_{\omega}\cdot\|_{k}$ to $\tilde{\mathcal{A}}$, respectively, and then utilizing \eqref{202109221247nn}$_1$, product estimate \eqref{product21}, and \eqref{small1}, we  arrive at \eqref{al1}--\eqref{alvarpi1}.
\hfill$\Box$
\end{pf}

The following lemma gives bounds on the divergence of $(u,\eta,\partial_{\omega}\eta)$.
\begin{lem}\label{lem:1400n}
Under the assumption \eqref{small1} 
with sufficiently small $\delta$, there holds that
\begin{align}
&\label{202209210921}
\|\mm{div}\eta\|_{i}\lesssim_{0}\|\eta\|_{3}\|\eta\|_{i+1}\;\;\hbox{ for }\;0\leqslant i\leqslant12,\\[1mm]
&\label{202209210921n}
\|\mm{div}u\|_{i}\lesssim_{0}\|\eta\|_{3}\|u\|_{i+1}+o(i-2)\|\eta\|_{i+1}\|u\|_{3}\;\;\hbox{ for }\;0\leqslant i\leqslant11,\\[1mm]
&\label{202209210921nn}
\|\mm{div}\partial_{\omega}\eta\|_{i}\lesssim_{0}\|\eta\|_{3}\|\partial_{\omega}\eta\|_{i+1}
+o(i-2)\|\eta\|_{i+1}\|\partial_{\omega}\eta\|_{3}\;\;\hbox{ for }\;0\leqslant i\leqslant11,
\end{align}
where
$$o(s):=\begin{cases}
1
& \;\hbox{ for }\; 0\leqslant s;\\
0
& \;\hbox{ for }\; s<0.
\end{cases}$$
\end{lem}
\begin{pf}
We can refer to \cite[Lemma 3.2]{JFJSO2021} for the proof.
Now we provide a brief sketch.
By \eqref{202207011433n}, we make use of determinant expansion theorem to see that
\begin{align*}
1=\mm{det}\nabla\zeta=\mm{det}(\nabla\eta+I)
=1+\mm{div}\eta+r_2^{\eta}+r_3^{\eta},
\end{align*}
where
$$\begin{aligned}
&r_2^{\eta}:=-\partial_3\eta_2\partial_2\eta_3-\partial_3\eta_1\partial_1\eta_3-\partial_2\eta_1\partial_1\eta_2
+\partial_2\eta_2\partial_3\eta_3+\partial_1\eta_1\partial_3\eta_3+\partial_1\eta_1\partial_2\eta_2,\\[1mm]
&r_3^{\eta}:=\partial_1\eta_1(\partial_2\eta_2\partial_3\eta_3-\partial_2\eta_3\partial_3\eta_2)
-\partial_2\eta_1(\partial_1\eta_2\partial_3\eta_3-\partial_1\eta_3\partial_3\eta_2)
+\partial_3\eta_1(\partial_1\eta_2\partial_2\eta_3-\partial_1\eta_3\partial_2\eta_2).
\end{aligned}$$
Consequently,
\begin{align}
&\mm{div}\eta=-(r_2^{\eta}+r_3^{\eta}),\nonumber\\[1mm]
&\mm{div}\partial_{\omega}\eta=-(\partial_{\omega}r_2^{\eta}+\partial_{\omega}r_3^{\eta}).\nonumber
\end{align}
Applying $\|\cdot\|_{i}$ to the above identities, and then utilizing product estimate \eqref{product21} and the assumption \eqref{prio1}, we immediately obtain \eqref{202209210921} and \eqref{202209210921nn}.

In addition, by \eqref{202109221247nn}$_3$, one has
$$\mm{div}u=-\mm{div}_{\tilde{\ml{A}}}u=-\tilde{\ml{A}}:\nabla u,$$
where $\mm{div}_{\tilde{\ml{A}}}u$ is defined as $\mm{div}_{\ml{A}}u$ with $\tilde{\ml{A}}$ in place of $\ml{A}$.
Applying $\|\cdot\|_{i}$ to the above identity, and then utilizing \eqref{product21} and \eqref{al1}, the estimate \eqref{202209210921n} follows. This completes the proof.
\hfill$\Box$
\end{pf}

Finally, we provide the pressure estimate.
\begin{lem}\label{lem2209141510}
Under the assumption \eqref{small1} with sufficiently small $\delta$, we have
\begin{align}\label{qestimate}
\|q\|_{k+2}\lesssim_{0}
\begin{cases}
\|u\|_{2}\|u\|_3+m^2\mathcal{F}_1
& \;\hbox{ for }\; k=1;\\[1mm]
\|u\|_{2}\|u\|_{k+2}+m^2\mathcal{F}_{k}+
\|\eta\|_{k+2}\left(\|u\|_{2}\|u\|_3+m^2\mathcal{F}_1\right)
& \;\hbox{ for }\; 2\leqslant k\leqslant10,
\end{cases}
\end{align}
where
$$
\mathcal{F}_{k}:=\begin{cases}
\|\partial_{\omega}\eta\|_{2}\|\partial_{\omega}\eta\|_3
+\|\eta\|_3\|\partial_{\omega}^2\eta\|_2\;\mbox{or}\;\|\eta\|_2\|\partial_{\omega}^2\eta\|_3
& \;\hbox{ for }\; k=1;\\[1mm]
\|\partial_{\omega}\eta\|_{3}\|\partial_{\omega}\eta\|_{k+1}+\|\eta\|_3\|\partial_{\omega}^2\eta\|_{k+1}
+\|\eta\|_{k+1}\big(\|\partial_{\omega}^2\eta\|_{3}+\|\partial_{\omega}\eta\|_{3}^2\big)
& \;\hbox{ for }\; 2\leqslant k\leqslant10.
\end{cases}$$
\end{lem}
\begin{pf}
Please refer to \cite[Lemma 3.3]{JFJSO2021} for the proof.
\hfill$\Box$
\end{pf}

\subsection{Energy estimate}

From now on, we assume that $(\eta,u)$ satisfies \eqref{prio1} and \eqref{prio2} with a sufficiently small $\delta\in(0,1]$, which yields that
\begin{align}
&\label{202209141614}
K_1:=\sup_{t\in[0,T]}\|\eta\|_{11}\lesssim_{0}\|\eta\|_{9}^{1/2}\|\eta\|_{13}^{1/2}\lesssim_{\omega}\delta\in(0,1],\\[1mm]
&\label{202209141615}
K_2:=\sqrt{\sup_{t\in[0,T]}(\|\eta\|_{9}+\|\partial_{\omega}\eta\|_{12})(\|u\|_{12}+\|\eta\|_{13})}\lesssim_{\omega}\delta. 
\end{align}
It should be noted that the smallness of $\delta$ only depends on the parameter $\mu$ and the unit vector $\omega$.
Moreover, \eqref{small1} automatically holds under the assumptions of \eqref{prio1} and \eqref{prio2}.

Now we proceed to derive the energy estimate for $(\eta,u)$.
Let integer $i$ satisfy $0\leqslant i\leqslant 8$, and let $\alpha$ satisfy $|\alpha|=i$.
Applying $\partial^{\alpha}$ to \eqref{202109221247nn}$_2$ then yields
\begin{align}\label{202209141812}
\partial^{\alpha}\left(u_t- \mu\Delta u-m^2\partial_{\omega}\eta\right)
=\partial^{\alpha}\left(\mathcal{N}^{u}-\nabla_{\tilde{\ml{A}}}q-\nabla q\right),
\end{align}
where we have defined that
\begin{align*}
&\mathcal{N}^{u}:=\partial_{l}\left(\mathcal{N}^{u}_{1l},\mathcal{N}_{2l},\mathcal{N}_{3l}\right)^{\mm{T}},\;\;
\mathcal{N}^{u}_{jl}:=\mu\left(\ml{A}_{kl}\tilde{\ml{A}}_{km}+\tilde{\ml{A}}_{ml}\right)\partial_{m}u_{j}\;\;\mbox{and}\;\;
\nabla_{\tilde{\ml{A}}}q:=\tilde{\mathcal{A}}_{ij}\partial_{j}q.
\end{align*}

\begin{lem}\label{lem2209141510l}
Under both the assumptions \eqref{prio1} and \eqref{prio2} with sufficiently small $\delta$, it holds that
\begin{align}
&\label{202209141623}
\frac{\mm{d}}{\mm{d}t}\|\nabla^{i}(u,m\partial_{\omega}\eta)\|_{0}^2+\|u\|_{i+1}^2 \lesssim\left(K_2+K_2^2+K^2/m\right)\mathcal{D}_{i},\quad0\leqslant i\leqslant8.
\end{align}
\end{lem}
\begin{pf}
Taking the inner product of \eqref{202209141812} and  $\partial^{\alpha}u$ in  $L^2$,
and then integrating by parts over $\Omega$,
we have
\begin{align}
\label{202109231702n}
&\frac{1}{2}\frac{\mm{d}}{\mm{d}t}\left(\|\partial^{\alpha} u\|^2_0+
\|m\partial_{\omega}\partial^{\alpha} \eta\|_0^2\right)
+\| \partial^{\alpha}\nabla u\|_0^2\nonumber\\
&=-\mu\int \partial^{\alpha}\mathcal{N}^{u}_{jl}\partial_{l}\partial^{\alpha} u_{j}\mm{d}y
+\int \left(\partial^{\alpha}{q}\partial^{\alpha}\mm{div} u
+\partial^{\alpha}(\tilde{\ml{A}}_{ij}q)\partial_{j} \partial^{\alpha}u_{i}
\right)\mm{d}y
=:H_{1}+H_{2}.
\end{align}
Exploiting H\"older's inequality, \eqref{product21} and \eqref{al1}, we obtain
\begin{align}\label{202209141845}
H_1&\lesssim\|\mathcal{N}^{u}_{jl}\|_{i}\|\partial_{l} u\|_{i}\nonumber\\[1mm]
&\lesssim\|u\|_{i+1}\big(1+\|\tilde{\ml{A}}\|_2\big)\big(\|\tilde{\ml{A}}\|_2\|u\|_{i+1}+o(i-2)\|\tilde{\ml{A}}\|_{i}\|u\|_{3}\big)
\nonumber\\[1mm]
&\lesssim o(i-2)\|\eta\|_{i+1}\|u\|_3\|u\|_{i+1}+\|\eta\|_{3}\|u\|_{i+1}^2\nonumber\\[1mm]
&\lesssim\delta\|u\|_{i+1}^2.
\end{align}
In the same manner,
\begin{align}\label{202209141846}
H_2&\lesssim\left(\|\tilde{\ml{A}}\|_2\|u\|_{i+1}+o(i-2)\|\tilde{\ml{A}}\|_{i}\|u\|_{3}\right)\|q\|_{i}
+\|\tilde{\ml{A}}\|_{i}\|u\|_{i+1}\|q\|_{3}
\nonumber\\[1mm]
&\lesssim\left(\|\eta\|_3\|u\|_{i+1}+o(i-2)\|\eta\|_{i+1}\|u\|_{3}\right)\|q\|_{i}
+\|\eta\|_{i+1}\|u\|_{i+1}\|q\|_{3}.
\end{align}
where we have used \eqref{202209210921n}. Next, we divide the estimate of $H_2$ into two cases.

\textbf{\emph{(1) Case $0\leqslant i\leqslant2$.}}

Making use of \eqref{poincareg} and the interpolation inequality, we deduce that
\begin{align}
&\label{202209161014n}
\|\partial_{\omega}\eta\|_6\|\partial_{\omega}\eta\|_3
\lesssim\|\partial_{\omega}\eta\|_{12-i}^{1/2}\|\partial_{\omega}\eta\|_{6-i}^{1/2}\|\partial_{\omega}\eta\|_{i}
\lesssim\|\partial_{\omega}\eta\|_{12-i}\|\partial_{\omega}\eta\|_{i}
\quad\;\hbox{ for }\; 0 \leqslant i\leqslant 2.
\end{align}
By virtue of \eqref{poincareg}, \eqref{qestimate} and \eqref{202209161014n}, it follows from \eqref{202209141846} that
\begin{align}\label{202209150904}
H_2&\lesssim\|\eta\|_3\|u\|_{i+1}\|q\|_{i}+\|\eta\|_{i+1}\|u\|_{i+1}\|q\|_{3}\lesssim\|\eta\|_{3}\|u\|_{i+1}\|q\|_{3}
\nonumber\\[1mm]
&\lesssim\|\eta\|_{3}\|u\|_{i+1}\big(\|u\|_{2}\|u\|_3+m^2\mathcal{F}_1\big)\nonumber\\[1mm]
&\lesssim\|\eta\|_{3}\|u\|_{i+1}\big(\|u\|_{2}\|u\|_3
+{m^2(\|\partial_{\omega}\eta\|_2\|\partial_{\omega}\eta\|_3+\|\eta\|_3\|\partial_{\omega}^2\eta\|_2)}\big)\nonumber\\[1mm]
&\lesssim\|\eta\|_{3}\|u\|_{i+1}\|u\|_{2}\|u\|_3
+m^2\|\eta\|_{3}\|u\|_{i+1}\|\partial_{\omega}\eta\|_3\|\partial_{\omega}\eta\|_6
\nonumber\\[1mm]&
\lesssim(K_2^2+K^2/m)\mathcal{D}_{i}.
\end{align}

\textbf{\emph{(2) Case $3\leqslant i\leqslant8$.}}

Thanks to \eqref{poincareg} and \eqref{qestimate}, it follows from \eqref{202209141846} that
\begin{align}\label{202209150934}
H_2&\lesssim\left(\|\eta\|_3\|u\|_{i+1}+\|\eta\|_{i+1}\|u\|_{3}\right)\|q\|_{i}+\|\eta\|_{i+1}\|u\|_{i+1}\|q\|_{3}\nonumber\\[1mm]
&\lesssim \left(\|\eta\|_3\|u\|_{i+1}+\|\eta\|_{i+1}\|u\|_{3}\right)\big(\|u\|_{2}\|u\|_{i}+m^2\mathcal{F}_{i-2}
+\|\eta\|_{i}\left(\|u\|_{2}\|u\|_3+m^2\mathcal{F}_1\right)
\big)\nonumber\\[1mm]
&\quad\quad+\|\eta\|_{i+1}\|u\|_{i+1}\big(\|u\|_{2}\|u\|_3+m^2\mathcal{F}_1\big)\nonumber\\[1mm]
&\lesssim\left(\|\eta\|_3\|u\|_{i+1}+\|\eta\|_{i+1}\|u\|_{3}\right)\big(\|u\|_{2}\|u\|_{i}+m^2\mathcal{F}_{i-2}\big)\nonumber\\[1mm]
&\quad\quad+\|\eta\|_{i+1}\|u\|_{i+1}\big(\|u\|_{2}\|u\|_3+m^2\mathcal{F}_1\big)\nonumber\\[1mm]
&\lesssim\|\eta\|_{i+1}\|u\|_{i+1}(1+\|\eta\|_{i})
\big(\|u\|_{2}\|u\|_{i+1}+m^2\|\partial_{\omega}\eta\|_6\|\partial_{\omega}\eta\|_{3}\big)\nonumber\\[1mm]
&\quad\;+m^2\|\eta\|_{i+1}\|u\|_{i+1}\big(\|\eta\|_3\|\partial_{\omega}\eta\|_{i}
+\|\partial_{\omega}\eta\|_{3}\|\partial_{\omega}\eta\|_{i-1}\big)\nonumber\\[1mm]
&\quad\;+m^2\|\eta\|_{i+1}\|u\|_{3}\|\eta\|_{i-1}\big(\|\partial_{\omega}^2\eta\|_{3}+\|\partial_{\omega}\eta\|_{3}^2\big)\nonumber\\[1mm]
&\lesssim(K_2^2+K^2/m)\mathcal{D}_{i}
+m^2\|\eta\|_{i+1}\|\eta\|_{i-1}\|u\|_{3}\|\partial_{\omega}\eta\|_{4}
\lesssim(K_2^2+K^2/m)\mathcal{D}_{i},
\end{align}
where we have used \eqref{202209161014} in the last inequality.

Consequently, inserting \eqref{202209141845}, \eqref{202209150904}--\eqref{202209150934} into \eqref{202109231702n} then \eqref{202209141623} follows immediately.
This completes the proof.
\hfill$\Box$
\end{pf}

\begin{lem}\label{lem2209141642}
Under both the assumptions \eqref{prio1} and \eqref{prio2} with sufficiently small $\delta$, it holds that
\begin{align}
\label{202209141642n}
&\frac{\mm{d}}{\mm{d}t}\bigg(\frac{\mu}{2}\|\nabla^{i+1}\eta\|_{0}^2+\sum_{|\alpha|=i}\int\partial^{\alpha}\eta\cdot\partial^{\alpha}u\mm{d}y\bigg)
+\|m\partial_{\omega}\eta\|_{i}^2 \nonumber\\
&\lesssim\left(K_2^2+K^2/m+K/m\right)\mathcal{D}_{i}+\|u\|_{i}^2,\quad0\leqslant i\leqslant8.
\end{align}
\end{lem}
\begin{pf}
Taking the inner product of \eqref{202209141812} by $\partial^{\alpha}\eta$ in  $L^2$,
and then integrating by parts over $\Omega$,
we have
\begin{align}
\label{202109231702nnnnn}
&\frac{1}{2}\frac{\mm{d}}{\mm{d}t}\left(\|\nabla\partial^{\alpha} \eta\|^2_0/2+
\int\partial^{\alpha}\eta\cdot\partial^{\alpha}u\mm{d}y\right)
+\|m\partial^{\alpha}\partial_{\omega}\eta\|_0^2\nonumber\\
&=\|\partial^{\alpha} u\|^2_0-\mu\int \partial^{\alpha}\mathcal{N}^{u}_{jl}\partial_{l}\partial^{\alpha} \eta_{j}\mm{d}y
+\int \left(\partial^{\alpha}{q}\mm{div}\partial^{\alpha} \eta
+\partial^{\alpha}(\tilde{\ml{A}}_{ij}q)\partial_{j} \partial^{\alpha}\eta_{i}\right)\mm{d}y\nonumber\\
&=:\|\partial^{\alpha} u\|^2_0+H_{3}+H_{4}.
\end{align}
Similarly to \eqref{202209141845}--\eqref{202209141846}, we utilize H\"older's inequality, \eqref{202209141614}, \eqref{product21}, \eqref{al1} and \eqref{202209210921} to bound
\begin{align}
&\label{202209141846N}
H_3\lesssim\|\eta\|_{3}\|\eta\|_{i+1}\|u\|_{i+1}+o(i-2)\|u\|_3\|\eta\|_{i+1}^2,\\[1mm]
&\label{202209141846Nn}
H_4\lesssim\|\eta\|_{i+1}\|\eta\|_{3}\|q\|_{i}+o(i-2)\|\eta\|_{i+1}^2\|q\|_{3}.
\end{align}
Moreover, by  \eqref{poincareg} and the interpolation inequality, we infer that
\begin{align}
\label{202209161005}
\|\eta\|_{i+1}\|\eta\|_{3}
&\lesssim
\begin{cases}
\|\partial_{\omega}\eta\|_{8+i}^{1/2}\|\partial_{\omega}\eta\|_{i}^{1/2}
\|\partial_{\omega}\eta\|_{12-i}^{1/2}\|\partial_{\omega}\eta\|_{i}^{1/2}
& \;\hbox{ for }\; 0\leqslant i\leqslant2;\\[1mm]
\|\eta\|_{5+i}^{1/2}\|\eta\|_{i-3}^{1/2}
\|\partial_{\omega}\eta\|_{12-i}^{1/2}\|\partial_{\omega}\eta\|_{i}^{1/2}
& \;\hbox{ for }\; 3 \leqslant i\leqslant 8
\end{cases}\nonumber\\[1mm]
&\lesssim
\begin{cases}
\|\partial_{\omega}\eta\|_{12}\|\partial_{\omega}\eta\|_{i}
& \;\hbox{ for }\; 0\leqslant i\leqslant2;\\[1mm]
\|\eta\|_{5+i}^{1/2}\|\partial_{\omega}\eta\|_{12-i}^{1/2}\|\partial_{\omega}\eta\|_{i}
& \;\hbox{ for }\; 3 \leqslant i\leqslant 8
\end{cases}
\end{align}
and
\begin{align}
\label{202209161014}
&\|\eta\|_{i+1}^2
\lesssim\|\eta\|_{5+i}\|\eta\|_{i-3}
\lesssim\|\eta\|_{5+i}\|\partial_{\omega}\eta\|_{i}
\quad\quad\quad\quad\;\;\hbox{ for }\; 3 \leqslant i\leqslant 8.
\end{align}
Hence, it follows from \eqref{202209141846N}, and \eqref{202209161005}--\eqref{202209161014} that
\begin{align}\label{202209151824}
H_3 
&
\lesssim
\begin{cases}
\|\eta\|_{3}\|\eta\|_{i+1}\|u\|_{i+1}\lesssim\|\partial_{\omega}\eta\|_{12}\|\partial_{\omega}\eta\|_{i}\|u\|_{i+1}
& \;\hbox{ for }\; 0\leqslant i\leqslant2;\\[1mm]
\|\eta\|_{i+1}^2\|u\|_{i+1}\lesssim\|\eta\|_{5+i}\|\partial_{\omega}\eta\|_{i}\|u\|_{i+1}
& \;\hbox{ for }\; 3 \leqslant i\leqslant 8
\end{cases}\nonumber\\[1mm]
&\lesssim
K/m\mathcal{D}_{i}
\qquad\qquad\qquad\qquad\qquad\qquad\qquad\qquad\;\;\;\hbox{ for }\; 0 \leqslant i\leqslant 8.
\end{align}

Similarly to the derivation of the estimate of $H_2$, we divide the estimate of $H_{4}$ into two cases as well.

\textbf{\emph{(1) Case $0\leqslant i\leqslant2$.}}

Thanks to \eqref{poincareg}, \eqref{qestimate}, \eqref{202209161014n}, \eqref{202209161005} and the interpolation inequality, we deduce from \eqref{202209141846Nn} that
\begin{align}\label{202209151946}
H_4&\lesssim\|\eta\|_3\|\eta\|_{i+1}\|q\|_{i}+\|\eta\|_{i+1}^2\|q\|_{3}\lesssim\|\eta\|_3\|\eta\|_{i+1}\|q\|_{3}\nonumber\\[1mm]
&\lesssim\|\eta\|_{3}\|\eta\|_{i+1}\big(\|u\|_{2}\|u\|_3+m^2\mathcal{F}_1\big)\nonumber\\[1mm]
&\lesssim\|\eta\|_{3}\|\eta\|_{i+1}\big(\|u\|_{2}\|u\|_3
+{m^2(\|\partial_{\omega}\eta\|_2\|\partial_{\omega}\eta\|_3+\|\eta\|_3\|\partial_{\omega}^2\eta\|_2)}\big)\nonumber\\[1mm]
&\lesssim\|\partial_{\omega}\eta\|_{12}\|\partial_{\omega}\eta\|_{i}
\big(\|u\|_{5-i}\|u\|_{i+1}+\|m\partial_{\omega}\eta\|_{12-i}\|m\partial_{\omega}\eta\|_{i}\big)
\nonumber\\[1mm]
&\lesssim K^2/m\mathcal{D}_{i}.
\end{align}

\textbf{\emph{(2) Case $3\leqslant i\leqslant8$.}}

In the light of \eqref{poincareg}, \eqref{qestimate}, \eqref{202209161014n}, \eqref{202209161005} and \eqref{202209161014}, we obtain  that
\begin{align}\label{202209160956}
H_4&\lesssim\|\eta\|_{i+1}\|\eta\|_{3}\|q\|_{i}+\|\eta\|_{i+1}^2\|q\|_{3}\nonumber\\[1mm]
&\lesssim\|\eta\|_{i+1}\|\eta\|_{3}\big(\|u\|_{2}\|u\|_{i}+m^2\mathcal{F}_{i-2}
\big)
+\|\eta\|_{i+1}^2\big(\|u\|_{2}\|u\|_3+m^2\mathcal{F}_1\big)\nonumber\\[1mm]
&\lesssim\|\eta\|_{5+i}^{1/2}\|\partial_{\omega}\eta\|_{12-i}^{1/2}
\|\partial_{\omega}\eta\|_{i}\big(\|u\|_{2}\|u\|_{i}+m^2\mathcal{F}_{i-2}\big)
+\|\eta\|_{5+i}\|\partial_{\omega}\eta\|_{i}\big(\|u\|_{2}\|u\|_3+m^2\mathcal{F}_1\big)\nonumber\\[1mm]
&\lesssim(\|\eta\|_{5+i}+\|\partial_{\omega}\eta\|_{10})\|\partial_{\omega}\eta\|_{i}\big(\|u\|_{2}\|u\|_{i}
+m^2\|\partial_{\omega}\eta\|_{3}
\|\partial_{\omega}\eta\|_{6}\big)\nonumber\\[1mm]
&\quad+m^2\|\eta\|_{5+i}^{1/2}\|\partial_{\omega}\eta\|_{12-i}^{1/2}\|\partial_{\omega}\eta\|_{i}
\big(\|\partial_{\omega}\eta\|_{3}\|\partial_{\omega}\eta\|_{i-1}+\|\eta\|_3\|\partial_{\omega}\eta\|_{i}
+\|\eta\|_{i-1}\big(\|\partial_{\omega}^2\eta\|_{3}+\|\partial_{\omega}\eta\|_{3}^2\big)\big)
\nonumber\\[1mm]
&\lesssim(K_2^2+K^2/m)\mathcal{D}_{i}
+m^2\|\eta\|_{5+i}^{1/2}\|\partial_{\omega}\eta\|_{12-i}^{1/2}\|\partial_{\omega}\eta\|_{i}
(\|\eta\|_{i-1}\|\partial_{\omega}^2\eta\|_{3})\nonumber\\[1mm]
&\lesssim(K_2^2+K^2/m)\mathcal{D}_{i}
+m^2\|\eta\|_{5+i}^{3/2}\|\partial_{\omega}\eta\|_{12-i}^{1/2}\|\partial_{\omega}\eta\|_{i}^2
\nonumber\\[1mm]&
\lesssim(K_2^2+K^2/m)\mathcal{D}_{i}.
\end{align}
Finally, plugging \eqref{202209151824}--\eqref{202209160956} into \eqref{202109231702nnnnn} then immediately yields
\eqref{202209141642n}. This completes the proof of Lemma \ref{lem2209141642}.
\hfill$\Box$
\end{pf}

\subsection{Highest-order derivative energy evolution}
{
Differing from the highest-order energy estimates in \cite{JFJSO2021}, where the highest-order derivatives of the velocity increased with the magnetic field intensity $m$, the $m$-increasing terms are concealed in the integral term
 $\int\partial^{\alpha}\nabla q\cdot\partial^{\alpha}u\mm{d}y=
-\int\partial^{\alpha} q\partial^{\alpha}\mm{div}u\mm{d}y$.
To address this $m$-increasing case, we instead to utilize the vorticity estimate method in the subsequent analysis.
}

Applying $\partial^{\alpha}\mm{curl}$ with $|\alpha|=11$ to \eqref{202109221247nn}$_2$ yields
\begin{align}\label{202210051356}
\partial^{\alpha}\mm{curl}\left(u_t- \mu\Delta u-m^2\partial_{\omega}^2\eta\right)
=\partial^{\alpha}\mm{curl}\left(\mathcal{N}^{u}-\nabla_{\tilde{\ml{A}}}q\right).
\end{align}
The following lemma gives bounds on the vorticity of $(u,\eta,\partial_{\omega}\eta)$.
\begin{lem}\label{lem:1402}
Under both the assumptions \eqref{prio1} and \eqref{prio2} with sufficiently small $\delta$, we have
\begin{align}
&\label{202210051352nn}
\frac{1}{2}\frac{\mm{d}}{\mm{d}t}\|\big(\mm{curl}u,m\mm{curl}\partial_{\omega}\eta\big)\|_{11}^2
+\mu\|\mm{curl}u\|_{12}^2\lesssim\sqrt{\mathcal{E}_{H}\mathcal{E}_{3}(1+\mathcal{D}_{6})\mathcal{D}_{H}}
+\delta(\|u\|_{13}^2+\|m\partial_{\omega}\eta\|_{12}^2),\\[1mm]
&\label{202210051352nnnn}
\sum_{|\alpha|=11}\frac{\mm{d}}{\mm{d}t}
\left(\int\partial^{\alpha}\mm{curl}\eta\cdot\partial^{\alpha}\mm{curl}u\mm{d}y
+\frac{\mu}{2}\|\nabla\partial^{\alpha}\mm{curl}\eta\|_0^2\right)
+\|m\mm{curl}\partial_{\omega}\eta\|_{11}^2\nonumber\\
&\lesssim\mathcal{E}_{H}\sqrt{\mathcal{D}_{3}}+\sqrt{\mathcal{E}_{H}\mathcal{E}_{3}\mathcal{D}_{H}}+\mathcal{E}_{H}\mathcal{D}_{6}+\|u\|_{12}^2.
\end{align}
\end{lem}
\begin{pf}
Taking the inner product of  \eqref{202210051356} and $\partial^{\alpha}\mm{curl}u$ in $L^2$,  one has
\begin{align}\label{202210051424}
&\frac{1}{2}\frac{\mm{d}}{\mm{d}t}\|\partial^{\alpha}(\mm{curl}u,m\mm{curl}\partial_{\omega}\eta)\|_{0}^2
+\mu\|\nabla\partial^{\alpha}\mm{curl}u\|_{0}^2\nonumber\\[1mm]
&=-\int\partial^{\alpha^{-}}\mm{curl}\mathcal{N}_{u}\cdot\partial^{\alpha^{+}}\mm{curl}u\mm{d}y
+\int\partial^{\alpha^{-}}\mm{curl}(\nabla_{\tilde{\ml{A}}}q)\cdot\partial^{\alpha^{+}}\mm{curl}u\mm{d}y
=:I_{1}+I_{2},
\end{align}
where for some nonzero component $\alpha_{j}$ of $\alpha$, we have denoted $\alpha^{-}$ and $\alpha^{+}$ by
\begin{align}\label{202210141454}
\alpha^{-}_{k}=\alpha^{+}_{k}=\alpha_{k}\quad\mbox{for}\quad k\neq j,\quad\quad
\alpha^{-}_{j}=\alpha_{j}-1\quad\mbox{and}\quad \alpha^{+}_{j}=\alpha_{j}+1.
\end{align}
By virtue of  H\"older's inequality, \eqref{product21} and \eqref{al1}, we can estimate that
$$
\begin{aligned}
I_1&\lesssim_{0}\|\mathcal{N}^{u}\|_{11}\|u\|_{13} \lesssim_{\mu}\|\tilde{\ml{A}}(\tilde{\ml{A}}+I)\nabla u\|_{12}\|u\|_{13}
\\[1mm]
&\lesssim\left(\|\eta\|_3\| u\|_{13}+\|\eta\|_{13}\|u\|_3\right)\|u\|_{13}
\lesssim\sqrt{\mathcal{E}_{H}\mathcal{E}_{3}\mathcal{D}_{H}}+\delta\|u\|_{13}^2
\end{aligned}
$$
and
$$
\begin{aligned}
I_2&\lesssim\|\nabla_{\tilde{\ml{A}}}q\|_{11}\|u\|_{13}
\lesssim\left(\|\eta\|_3\| q\|_{12}+\|\eta\|_{12}\| q\|_{3}\right)\|u\|_{13}\\[1mm]
&\lesssim\|\eta\|_3\|u\|_{13}\big(\|u\|_{2}\|u\|_{12}+m^2\mathcal{F}_{10}\big)
+\|\eta\|_{12}\|u\|_{13}(1+\|\eta\|_3)\big(\|u\|_{2}\|u\|_3+m^2\mathcal{F}_1\big)\nonumber\\[1mm]
&\lesssim\|\eta\|_3\sqrt{\mathcal{E}_{H}\mathcal{E}_{3}\mathcal{D}_{H}}
+(\|\eta\|_3\|m\partial_{\omega}\eta\|_{6})\|u\|_{13}\|m\partial_{\omega}\eta\|_{12}\\[1mm]
&\quad\;
+\|\eta\|_{11}\|m\partial_{\omega}\eta\|_{7}\|m\partial_{\omega}\eta\|_{3}\|u\|_{13}
+\|\eta\|_{11}\|\eta\|_{3}\|m\partial_{\omega}\eta\|_{3}^2\|u\|_{13}\\[1mm]
&\quad\;+\|\eta\|_{12}\|u\|_{13}\big(\|u\|_{2}\|u\|_3+\|m\partial_{\omega}\eta\|_{3}\|m\partial_{\omega}\eta\|_{6}\big)
\nonumber\\[1mm]
&\lesssim \sqrt{\mathcal{E}_{H}\mathcal{E}_{3}\mathcal{D}_{H}}
+\sqrt{\mathcal{E}_{H}\mathcal{D}_{H}\mathcal{E}_{3}\mathcal{D}_{6}}
+\delta(\|u\|_{13}^2+\|m\partial_{\omega}\eta\|_{12}^2),\nonumber\\[1mm]
\end{aligned}
$$
where we have used \eqref{202209141614} and the fact
$$
\begin{aligned}
&m^2\|\eta\|_{3}\|\partial_{\omega}\eta\|_{4}\lesssim\|m\partial_{\omega}\eta\|_{6}\|m\partial_{\omega}\eta\|_{4}
\lesssim \|m\partial_{\omega}\eta\|_{7}\|m\partial_{\omega}\eta\|_{3}.
\end{aligned}$$
Plugging the above two inequalities into \eqref{202210051424}, then \eqref{202210051352nn} follows immediately.

Similarly to the derivation of \eqref{202210051424},
we take the inner product of \eqref{202210051356} and $\partial^{\alpha}\mm{curl}\eta$ in $L^2$ to have
\begin{align}\label{202210051426}
&\frac{\mm{d}}{\mm{d}t}
\left(\int\partial^{\alpha}\mm{curl}\eta\cdot\partial^{\alpha}\mm{curl}u\mm{d}y
+\frac{\mu}{2}\|\nabla\partial^{\alpha}\mm{curl}\eta\|_0^2\right)
+\|m\partial^{\alpha}\mm{curl}\partial_{\omega}\eta\|_{0}^2\nonumber\\[1mm]
&=\|\partial^{\alpha}\mm{curl}u\|_0^2
-\int\partial^{\alpha^{-}}\mm{curl}\mathcal{N}^{u}\cdot\partial^{\alpha^{+}}\mm{curl}\eta\mm{d}y
+\int\partial^{\alpha^{-}}\mm{curl}(\nabla_{\tilde{\ml{A}}}q)\cdot\partial^{\alpha^{+}}\mm{curl}\eta\mm{d}y\nonumber\\[1mm]
&
=:\|\partial^{\alpha}\mm{curl}u\|_0^2+I_{3}+I_{4}.
\end{align}
Applying H\"older's inequality, \eqref{product21}, \eqref{al1} and Young's inequality then gives rise to
$$
\begin{aligned}
I_3&\lesssim_{0}\|\mathcal{N}^{u}\|_{11}\|\eta\|_{13} \lesssim_{\mu}\|\tilde{\ml{A}}(\tilde{\ml{A}}+I)\nabla u\|_{12}\|\eta\|_{13}
\\[1mm]
&\lesssim\left(\|\eta\|_3\| u\|_{13}+\|\eta\|_{13}\|u\|_3\right)\|\eta\|_{13}
\lesssim \mathcal{E}_{H}\sqrt{\mathcal{D}_{3}}+\sqrt{\mathcal{E}_{H}\mathcal{E}_{3}}\mathcal{D}_{H}
\end{aligned}
$$
and
$$
\begin{aligned}
I_4&\lesssim\|\nabla_{\tilde{\ml{A}}}q\|_{11}\|\eta\|_{13}
\lesssim\left(\|\eta\|_3\| q\|_{12}+\|\eta\|_{12}\| q\|_{3}\right)\|\eta\|_{13}\\[1mm]
&\lesssim\|\eta\|_3\|\eta\|_{13}\big(\|u\|_{2}\|u\|_{12}+m^2\mathcal{F}_{10}\big)
+(1+\|\eta\|_3)\|\eta\|_{12}\|\eta\|_{13}\big(\|u\|_{2}\|u\|_3+m^2\mathcal{F}_1\big)\nonumber\\[1mm]
&\lesssim\|\eta\|_3\mathcal{E}_{H}\sqrt{\mathcal{D}_{3}}
+\mathcal{E}_{H}\mathcal{E}_{3}/\mathcal{D}_{6}\\[1mm]
&\quad
+\|\eta\|_3\|\eta\|_{13}(m\|\eta\|_{3}\|m\partial_{\omega}\eta\|_{12}
+m\|\eta\|_{11}\|m\partial_{\omega}\eta\|_{4}+\|\eta\|_{11}\|m\partial_{\omega}\eta\|_{3}^2)\\[1mm]
&\quad+
\|\eta\|_{12}\big(\|u\|_{2}\|u\|_3+\|m\partial_{\omega}\eta\|_{3}\|m\eta\|_{3}\big)\|\eta\|_{13}
\nonumber\\[1mm]
&\lesssim
 \mathcal{E}_{H}\sqrt{\mathcal{D}_{3}}+\mathcal{E}_{H}\mathcal{D}_{6},\nonumber\\[1mm]
\end{aligned}
$$
where we have used \eqref{202209141614} and
$$
\begin{aligned}
&m^2\|\eta\|_{3}\|\partial_{\omega}\eta\|_{4}\lesssim\|m\partial_{\omega}\eta\|_{6}\|m\partial_{\omega}\eta\|_{4}
\lesssim\|m\partial_{\omega}\eta\|_{7}\|m\partial_{\omega}\eta\|_{3},\\
&m\|\eta\|_{3}^2\lesssim m\|\eta\|_{6}\|\eta\|_{0}\lesssim\|\eta\|_{6}\|m\partial_{\omega}\eta\|_{3}.
\end{aligned}$$
Putting both the two estimates of $I_3$ and $I_4$ into \eqref{202210051426} then yields \eqref{202210051352nnnn} immediately.
This completes the proof of Lemma \ref{lem:1402}.
\hfill$\Box$
\end{pf}
\subsection{Stability estimate}\label{stabilitye}
With Lemmas \ref{lem2209141510l}--\ref{lem:1402} 
 and Lemma \ref{lem:1400n} in  hand, we are now able to establish the a \emph{priori} estimate \eqref{pr} under the assumptions \eqref{prio1} and \eqref{prio2} with sufficiently small $\delta$.

First of all, it follows from \eqref{202209141623} and \eqref{202209141642n} with $0\leqslant i\leqslant8$ that
\begin{align*}
\frac{\mm{d}}{\mm{d}t}\tilde{\mathcal{E}}_{i}+c\mathcal{D}_{i}\lesssim\left(K_2+K_2^2+K/m+K^2/m\right)\mathcal{D}_{i},
\end{align*}
which implies by \eqref{prio2} and \eqref{202209141615} that
\begin{align}\label{202209171101}
\frac{\mm{d}}{\mm{d}t}\tilde{\mathcal{E}}_{i}+c\mathcal{D}_{i}\lesssim0\quad\;\mbox{for}\;0\leqslant i\leqslant8,
\end{align}
where we have defined that
$$\tilde{\mathcal{E}}_{i}=\bigg(\frac{\mu}{2}\|\nabla^{i+1} \eta\|^2_0+
\sum_{|\alpha|=i}\int\partial^{\alpha}\eta\cdot\partial^{\alpha}u\mm{d}y\bigg)
+C_1\|\nabla^{i}(u,m\partial_{\omega}\eta)\|_{0}^2.$$
Choosing a suitable positive constant $C_1$ above, and utilizing \eqref{poincareg}, one can see that
\begin{align}\label{202209171106}
\tilde{\mathcal{E}}_{i}\;\mbox{ is equivalent to }\;\mathcal{E}_{i}
\end{align}
and the following inequalities
\begin{align}\label{202209171114}
\tilde{\mathcal{E}}_{i}\lesssim\mathcal{E}_{i}\lesssim\langle m^{-1}\rangle^{2}\mathcal{D}_{i+4}\quad\mbox{for}\;\;0\leqslant i\leqslant 8.
\end{align}
Hence, we can further deduce from \eqref{202209171101} and \eqref{202209171114} that
\begin{align}\label{202209171117}
\sum_{j=0}^2\left(d_{j}\langle m^{-1}\rangle^{2j}\frac{\mm{d}}{\mm{d}t}\big(\langle t\rangle^{(2-j)}\tilde{\mathcal{E}}_{4j}\big)
+h_{j}\langle m^{-1}\rangle^{2j}\langle t\rangle^{(2-j)}\mathcal{D}_{4j}\right)\lesssim0
\end{align}
for some constants $d_{j}$ and $h_{j}$ depending on $\mu$ and $\omega$.
Consequently, integrating \eqref{202209171101} and \eqref{202209171117} in time over $(0,t)$ respectively, and then utilizing \eqref{202209171114}, we obtain
\begin{align}
&\label{202209171101n}
{\mathcal{E}}_{i}(t)+c\int_0^{t}\mathcal{D}_{i}(\tau)\mm{d}\tau\lesssim{\mathcal{E}}_{i}(0),\;\;0\leqslant i\leqslant8,\\[1mm]
&\label{202209171117n}
\sum_{j=0}^2\left(\langle m^{-1}\rangle^{2j}\big(\langle t\rangle^{(2-j)}{\mathcal{E}}_{4j}(t)\big)
+\langle m^{-1}\rangle^{2j}\int_{0}^{t}\langle \tau\rangle^{(2-j)}\mathcal{D}_{4j}(\tau)\mm{d}\tau\right)
\lesssim\Xi.
\end{align}

In addition, by exploiting the interpolation inequality along with  \eqref{202209171101n}--\eqref{202209171117n}, one has
\begin{align}\label{202209171135}
\mathcal{E}_{3}(t)\lesssim\mathcal{E}_{8}(t)^{3/8}\mathcal{E}_{0}(t)^{5/8}\lesssim\mathcal{E}_{8}(0)^{3/8}(\langle t\rangle^{-2}{\mathcal{E}}_{0}(0))^{5/8}\lesssim\Xi\langle t\rangle^{-5/4}.
\end{align}
Since
$$\begin{aligned}
&\frac{\mm{d}}{\mm{d}t}\left(\langle t\rangle^{\frac{5/4+1}{2}}\tilde{\mathcal{E}}_{3}\right)
+\langle t\rangle^{\frac{5/4+1}{2}}{\mathcal{D}}_{3}\\[1mm]
&=\langle t\rangle^{\frac{5/4+1}{2}}\frac{\mm{d}}{\mm{d}t}\tilde{\mathcal{E}}_{3}
+ \langle t\rangle^{\frac{5/4+1}{2}}{\mathcal{D}}_{3}+\frac{5/4+1}{2}\tilde{\mathcal{E}}_{3}\langle t\rangle^{\frac{1}{8}},
\end{aligned}$$
from \eqref{202209171101} with $i=3$, we deduce
\begin{align*}
&\frac{\mm{d}}{\mm{d}t}\left(\langle t\rangle^{\frac{5/4+1}{2}}\tilde{\mathcal{E}}_{3}\right)
+\langle t\rangle^{\frac{5/4+1}{2}}\tilde{\mathcal{D}}_{3}
\leqslant\frac{5/4+1}{2}\tilde{\mathcal{E}}_{3}\langle t\rangle^{\frac{1}{8}},
\end{align*}
which together with \eqref{202209171135} yields that
\begin{align*}
&\frac{\mm{d}}{\mm{d}t}\left(\langle t\rangle^{\frac{5/4+1}{2}}\tilde{\mathcal{E}}_{3}\right)
+\langle t\rangle^{\frac{5/4+1}{2}}{\mathcal{D}}_{3}
\lesssim\Xi\langle t\rangle^{-\frac{9}{8}}.
\end{align*}
Hence we arrive at
\begin{align}
\label{202109261645}
\sup_{0\leqslant \tau< t}\langle \tau\rangle^{\frac{5/4+1}{2}}\mathcal{E}_{3}(\tau)+\int_0^t{\langle \tau\rangle^{\frac{5/4+1}{2}}}{\mathcal{D}_{3}}(\tau)\mm{d}\tau
\lesssim({\mathcal{E}}_{3}(0)+\Xi)\lesssim\Xi.
\end{align}

Now we are in position to establish the subtle higher-order energy estimate. Similarly to \eqref{202209171101},
we deduce from \eqref{202209141623}, \eqref{202209141642n}, \eqref{202209210921}--\eqref{202209210921nn} and
Hodge type elliptic estimate \eqref{hodgeellipticnn} that
\begin{align*}
\frac{\mm{d}}{\mm{d}t}\tilde{\mathcal{E}}_{H}+c\mathcal{D}_{H}
&\lesssim
\sqrt{\mathcal{E}_{H}\mathcal{E}_{3}(1+\mathcal{D}_{6})\mathcal{D}_{H}}
+\mathcal{E}_{H}\sqrt{\mathcal{D}_{3}}+\mathcal{E}_{H}\mathcal{D}_{6},
\end{align*}
which gives rise to
\begin{align}\label{202209171626}
\frac{\mm{d}}{\mm{d}t}\tilde{\mathcal{E}}_{H}+c\mathcal{D}_{H}
\lesssim\mathcal{E}_{H}\big(\sqrt{\mathcal{D}_{3}}+\mathcal{E}_3+\left(1+\mathcal{E}_3\right)\mathcal{D}_{6}\big),
\end{align}
where we have defined that
$$\tilde{\mathcal{E}}_{H}:=\sum_{|\alpha|=11}\left(\int\partial^{\alpha}\mm{curl}\eta\cdot\partial^{\alpha}\mm{curl}u\mm{d}y
+\frac{\mu}{2}\|\nabla\partial^{\alpha}\mm{curl}\eta\|_0^2\right)
+C_2\|\big(\mm{curl}u,m\mm{curl}\partial_{\omega}\eta\big)\|_{11}^2.$$
By \eqref{202209210921}--\eqref{202209210921nn}, and
Hodge type elliptic estimate \eqref{hodgeellipticnn}, one sees that
\begin{align}\label{202209171106nn}
\tilde{\mathcal{E}}_{H}\;\;\mbox{is equivalent to}\;\;\mathcal{E}_{H}
\end{align}
for some suitable positive constant $C_2$.

With \eqref{202209171626}, \eqref{202109261645} and \eqref{202209171101n} in hand, we utilize Gronwall's argument to \eqref{202209171626} and then utilize \eqref{202209171101n} and the integrable decay of the lower dissipation \eqref{202109261645} to deduce that,
there exists a constant $\delta_1\in(0,1]$ such that, for any $\delta\leqslant\delta_1$,
\begin{align}\label{202209171712}
\mathcal{E}_{H}(t)
&\lesssim\mathcal{E}_{H}^{0}e^{c\int_0^{t}(\sqrt{\mathcal{D}_{3}}+\mathcal{E}_3+(1+\mathcal{E}_3)\mathcal{D}_{6})(\tau)\mm{d}\tau}\nonumber\\[1mm]
&\lesssim\mathcal{E}_{H}^{0}e^{c\left(\int_0^{t}\sqrt{\mathcal{D}_{3}}(\tau)\mm{d}\tau
+\int_0^{t}\mathcal{E}_3(\tau)\mm{d}\tau
+(1+\mathcal{E}_3(0))\int_0^{t}\mathcal{D}_{6}(\tau)\mm{d}\tau\right)}\nonumber\\[1mm]
&\lesssim \mathcal{E}_{H}^{0}e^{c\left(\Xi^{1/2}
+\Xi+(1+\mathcal{E}_{3}(0))\mathcal{E}_6(0)\right)}\nonumber\\[1mm]
&\leqslant c_1\mathcal{E}_{H}^{0}e^{c_2\left(\Xi^{1/2}
+\Xi+\mathcal{E}_{3}(0)\mathcal{E}_{6}(0)\right)}/4:=c_1\mathcal{E}_{H}^{0}e^{c_2\Psi}/4,
\end{align}
where $c_1\geqslant4$ and
\begin{align}\label{202210141600}
\Psi:=\Xi^{1/2}+\Xi+\mathcal{E}_{3}^0\mathcal{E}_{6}^0.
\end{align}
Furthermore, we deduce from \eqref{202209171626} and \eqref{202209171712} that
\begin{align}
&\label{202209171723}
{\mathcal{E}}_{H}(t)+\int_0^{t}\mathcal{D}_{H}(\tau)\mm{d}\tau\lesssim\mathcal{E}_{H}^0(1+e^{c_2\Psi}\Psi).
\end{align}

Now we let
\begin{align}\label{202209171746}
K:=\sqrt{c_1\mathcal{E}_{H}^{0}e^{c_2\Psi}}>0,
\end{align}
then we immediately obtain the desired a \emph{priori} estimate \eqref{pr} from \eqref{202209171712} under the a \emph{priori} assumptions \eqref{prio1}--\eqref{prio2} with $\delta\leqslant\delta_1$.

Consequently, with the \emph{priori} estimates \eqref{202209171101n}--\eqref{202209171117n}, \eqref{202109261645}, \eqref{202209171712} and \eqref{202209171723} (under the conditions \eqref{prio1}--\eqref{prio2} with $\delta\leqslant\delta_1$) in hand, and then combining with the local well-posedness result and a result concerning diffeomorphism mapping below,
we can  utilize a continuity argument as in \cite{JFJSO2021} to establish Theorem \ref{thm:2022n}.

\begin{pro}\label{0301lwp}
Let $(\eta^0,u^0)\in H^{13}\times H^{12}$ satisfy $\|(\nabla\eta^0,u^0)\|_{12}\leqslant B$ and $\mm{div}_{\ml{A}^0}u^0=0$, where $B$ is a positive constant, and $\ml{A}^0:=(\nabla\eta^0+I)^{-\mm{T}}$.
Then, there is a constant $\delta_2\in(0,1)$, such that for any $(\eta^0,u^0)$ satisfying
\begin{align}\label{202209291536}
\|\nabla\eta^0\|_{7}\leqslant\delta_2,
\end{align}
there exists a local time $T>0$ (depending possibly on $B$, $\mu$, $m$ and $\delta$) and a unique local classical solution
$(\eta,u,q)\in H^{13}\times H^{12}\times \underline{H}^{12}$ to \eqref{202109221247nn}, satisfying
$$0<\inf_{(y,t)\in\mathbb{R}^3\times[0,T]}\mm{det}(\nabla\eta+I)\quad
\mbox{and}\quad \sup_{t\in[0,T]}\|\nabla\eta\|_{7}\leqslant 2\delta_2.$$
\end{pro}
\begin{pf}
Please refer to \cite[Proposition 3.1]{JFJSO2020}.
\hfill$\Box$
\end{pf}

\begin{pro}\label{1014dm}
There exists a positive constant $\delta_3$ such that, for any $\varphi\in H^{13}$ satisfying $\|\nabla\varphi\|_{2}\leqslant\delta_3$, we have (after possibly being redefined on a set of zero measure) $\mm{det}(\nabla\varphi+I)>1/2$ and
\begin{align}
\psi:=\varphi+y:\mathbb{R}^3\rightarrow\mathbb{R}^3\;\mbox{is a}\;\;C^{11}\;\mbox{homeomorphism mapping}.
\end{align}
\end{pro}
\begin{pf}
Please refer to \cite[Lemma 4.2]{JFJSOMITIN}.
\hfill$\Box$
\end{pf}

\section{Proof of Theorem \ref{thm:2022nn}.
}\label{AB}
This section is devoted to the proof of Theorem \ref{thm:2022nn}.
Let us first recall $(\eta^{d}, u^{d})=(\eta,u)-(\eta^{L},u^{L})$, then $(\eta^{d}, u^{d})$ satisfies
\begin{equation}\label{20210922111nn}
\begin{cases}
\eta^{d}_t=u^{d}&\mbox{ in } \mathbb{T}^3 ,\\[1mm]
u^{d}_t - \mu\Delta u^{d}-m^2\partial_{\omega}^2\eta^{d}=\mathcal{N}^{u}-\nabla_{\mathcal{A}}q
&\mbox{ in } \mathbb{T}^3 ,\\[1mm]
\mm{div} u^{d}= -\mm{div}_{\tilde{\ml{A}}}u
&\mbox{ in } \mathbb{T}^3 ,\\[1mm]
(\eta^{d},u^{d})|_{t=0}=(\eta^{\mm{r}}, u^{\mm{r}})
&\mbox{ in } \mathbb{T}^3,
\end{cases}
\end{equation}
Moreover, from \eqref{20210922111nn} it follows that
$(\eta^{d})_{\mathbb{T}^3}=(u^{d})_{\mathbb{T}^3}=0$, and $\mm{div}\eta^{d}=\mm{div}\eta$ for any $t>0$,
since $(\eta^{\mm{r}})_{\mathbb{T}^3}=(u^{\mm{r}})_{\mathbb{T}^3}=0$ and $\mm{div}\eta^{L}=0$.
In what follows,
we define
$$\mathfrak{E}_{j}^{d}:=\|(\nabla\eta^{d}, u^{d}, m\partial_{\omega}\eta^{d})\|_{j}^2\quad\mbox{and}\;\;\;
\mathfrak{D}_{j}^{d}:=\|(\nabla u^{d}, m\partial_{\omega}\eta^{d})\|_{j}^2,\quad1\leqslant j\leqslant9.$$
It is easy to see that, for $j\geqslant2$,
\begin{align}\label{202210161627}
\begin{aligned}
&\limsup_{t\rightarrow0}\sup_{\tau\in[0,t]}\mathfrak{E}_{j}^{d}\\
&
\lesssim\|(\nabla\eta^{\mm{r}}, u^{\mm{r}}, m\partial_{\omega}\eta^{\mm{r}})\|_{j}^2\\[1mm]
&\lesssim\|\eta^0\|_{3}^2\|(\nabla\eta^{0}, u^{0}, m\partial_{\omega}\eta^{0})\|_{j}^2
+o(j-2)\|\eta^{0}\|_{j}^2\|(u^{0}, m\partial_{\omega}\eta^{0})\|_{3}^2\\[1mm]
&\lesssim m^{-2}\big(\|m\partial_{\omega}\eta^0\|_{6}^2\|(\nabla\eta^{0}, u^{0}, m\partial_{\omega}\eta^{0})\|_{j}^2
+o(j-2)\|m\partial_{\omega}\eta^{0}\|_{j+3}^2\|(u^{0}, m\partial_{\omega}\eta^{0})\|_{3}^2\big).
\end{aligned}
\end{align}

{
Keeping in mind that the authors in \cite{JFJSO2021} directly applied $\partial^{\alpha}$
to \eqref{20210922111nn}$_2$ to derive error energy estimates for $(\eta^{d}, u^{d})$, resulting in a subtle term
$\int\partial^{\alpha}\nabla q\cdot\partial^{\alpha}u^{d}\mm{d}y=-\int\partial^{\alpha}q\mm{div}\partial^{\alpha}u^{d}\mm{d}y$.
Indeed, in \cite{JFJSO2021}, the authors had to use \eqref{20210922111nn}$_3$ to replace $\mm{div}u^{d}$, which is given by
$$\int\partial^{\alpha}\nabla q\cdot\partial^{\alpha}u^{d}\mm{d}y=-\int\partial^{\alpha}q\mm{div}\partial^{\alpha}u^{d}\mm{d}y
=\int\partial^{\alpha}q\partial^{\alpha}\mm{div}_{\tilde{\ml{A}}}u\mm{d}y.$$
As a result, whether it is the higher-order error energy or the lower-order error energy, the convergence rate of the error energy is just $m^{-1}$. This is because the convergence rate in $m$ is later determined by $\|\partial^{\alpha}\mm{div}_{\tilde{\ml{A}}}u\|_{0}$ at this point. To enhance the convergence rate of the nonlinear system towards a linearized problem, we instead choose to use the $\mm{curl}$-estimate method as follows.
}

Following the same arguments as in the proof of Lemma \ref{lem:1402} by slight modifications,
we first apply $\partial^{\alpha}\mm{curl}$ to \eqref{20210922111nn}$_2$ to deduce that
\begin{align}\label{202210051721}
\partial^{\alpha}\mm{curl}(u^{d}_t - \mu\Delta u^{d}-m^2\partial_{\omega}^2\eta^{d})
=\partial^{\alpha}\mm{curl}(\mathcal{N}^{u}-\nabla_{\tilde{\mathcal{A}}}q),
\end{align}
where $\alpha$ satisfies $|\alpha|:=i$ with $i=0,4,7,8$.

Taking the inner product of \eqref{202210051721} by $\partial^{\alpha}\mm{curl}u^{d}$ and $\partial^{\alpha}\mm{curl}\eta^{d}$ in $L^2$, respectively, then integrating by parts over $\Omega$, and finally using \eqref{20210922111nn}$_1$, we obtain that
\begin{align*}
&\frac{1}{2}\frac{\mm{d}}{\mm{d}t}\|\partial^{\alpha}\big(\mm{curl}u^{d},m\partial_{\omega}\mm{curl}\eta^{d}\big)\|_{0}^2
+\mu\|\nabla\partial^{\alpha}\mm{curl}u^{d}\|_{0}^2\nonumber\\[1mm]
&=\int\partial^{\alpha}\mm{curl}\mathcal{N}^{u}\cdot\partial^{\alpha}\mm{curl}u^{d}\mm{d}y
-\int\partial^{\alpha}\mm{curl}(\nabla_{\tilde{\ml{A}}}q)\cdot\partial^{\alpha}\mm{curl}u^{d}\mm{d}y
=:J_{1}+J_{2}
\end{align*}
and
\begin{align*}
&\frac{\mm{d}}{\mm{d}t}
\left(\int\partial^{\alpha}\mm{curl}\eta^{d}\cdot\partial^{\alpha}\mm{curl}u^{d}\mm{d}y
+\frac{\mu}{2}\|\nabla\partial^{\alpha}\mm{curl}\eta^{d}\|_0^2\right)
+\|m\partial^{\alpha}\mm{curl}\partial_{\omega}\eta^{d}\|_{0}^2\nonumber\\[1mm]
&=\|\partial^{\alpha}\mm{curl}u^{d}\|_0^2
+\int\partial^{\alpha}\mm{curl}\mathcal{N}^{u}\cdot\partial^{\alpha}\mm{curl}\eta^{d}\mm{d}y
-\int\partial^{\alpha}\mm{curl}(\nabla_{\tilde{\ml{A}}}q)\cdot\partial^{\alpha}\mm{curl}\eta^{d}\mm{d}y\nonumber\\[1mm]
&
=:\|\partial^{\alpha}\mm{curl}u\|_0^2+J_3+J_4.
\end{align*}
From which it follows that
\begin{align}\label{202210151940}
\frac{\mm{d}}{\mm{d}t}\mathcal{E}^{d}_{\mm{curl},i}+\mathcal{D}^{d}_{\mm{curl},i}
\lesssim J_{1}+J_{2}+J_3+J_4,
\end{align}
where for some suitable large $C_3$ (depending on $\mu$),
\begin{align}
&\mathcal{E}^{d}_{\mm{curl},i}:=
\frac{\mu}{4}\|\nabla\partial^{\alpha}\mm{curl}\eta^{d}\|_0^2
+C_3\|\partial^{\alpha}(\mm{curl}u^{d},m\mm{curl}\partial_{\omega}\eta^{d})\|_{0}^2,\nonumber\\[1mm]
&\mathcal{D}^{d}_{\mm{curl},i}:=\|m\partial^{\alpha}\mm{curl}\partial_{\omega}\eta^{d}\|_{0}^2
+C_3\mu\|\nabla\partial^{\alpha}\mm{curl}u^{d}\|_{0}^2.\nonumber
\end{align}

Next, we estimate $J_1,\cdots, J_4$ in sequel,
and we shall split it into two cases, namely case $i=0$ and case $i\neq0$.
To begin with, by virtue of \eqref{202210061225}, the pressure estimate \eqref{qestimate} for $i\leqslant8$ can be simplified as follows.
\begin{align}\label{qestimate1}
\|q\|_{i+1}\lesssim_{0}
\begin{cases}
\|u\|_{2}\|u\|_{3}+\|m\partial_{\omega}\eta\|_4\|m\partial_{\omega}\eta\|_{5}
& \;\hbox{ for }\; i=0;\\[1mm]
\|u\|_{2}\|u\|_{i+1}+\|m\partial_{\omega}\eta\|_{6}\|m\partial_{\omega}\eta\|_{i+1}
+m\|\eta\|_{i}\|m\partial_{\omega}\eta\|_{4}
& \;\hbox{ for }\; i=4,7,8.
\end{cases}
\end{align}

\textbf{\emph{(1) Case $i=0$.}}

Integrating by part (or maybe twice) over $\Omega$, and then utilizing \eqref{poincareg}, \eqref{product21}, \eqref{al1}, \eqref{qestimate1} and the interpolation inequality, we infer that
\begin{align}
J_1 \lesssim&
\|\eta\|_3\|u\|_2\|\mm{curl}u^{d}\|_1
\lesssim m^{-1}\|m\partial_{\omega}\eta\|_0^{1/2}\|m\partial_{\omega}\eta\|_{12}^{1/2}\|u\|_0^{1/2}\|u\|_4^{1/2}\|\mm{curl}u^{d}\|_1 \nonumber\\[1mm]
\lesssim& C_{\epsilon}m^{-2}\langle t\rangle^{-2}(\langle t\rangle^{2}\mathcal{E}_{0})\mathcal{D}_{12}+\epsilon\|\mm{curl}u^{d}\|_1^2,\label{k11}\\[1mm]
J_2 \lesssim&
\|\eta\|_{3}\|\nabla q\|_0\|\mm{curl}u^{d}\|_1\nonumber\\[1mm]
\lesssim&\|\eta\|_{3}(\|u\|_{2}\|u\|_{3}+\|m\partial_{\omega}\eta\|_{4}\|m\partial_{\omega}\eta\|_{5})
\|\mm{curl}u^{d}\|_1\nonumber\\[1mm]
\lesssim& C_{\epsilon}m^{-2}\|m\partial_{\omega}\eta\|_{12}\|m\partial_{\omega}\eta\|_0
(\|u\|_{2}^2\|u\|_{3}^2
+\|m\partial_{\omega}\eta\|_{4}^2\|m\partial_{\omega}\eta\|_{5}^2)
+\epsilon\|\mm{curl}u^{d}\|_1^2\nonumber\\[1mm]
\lesssim& C_{\epsilon}m^{-2}\langle t\rangle^{-2}\mathcal{E}_{12}(\langle t\rangle^{2}\mathcal{E}_{0})\mathcal{D}_{12}+\epsilon\|\mm{curl}u^{d}\|_1^2.\label{k21}
\end{align}
In the same way,
\begin{align}
J_3 \lesssim&
\|\eta\|_3\|u\|_1\|\eta^{d}\|_3
\lesssim m^{-2}\|m\partial_{\omega}(\eta,\eta^{L})\|_{12}\underbrace{\|m\partial_{\omega}(\eta,\eta^{L})\|_0\|u\|_1}\nonumber\\[1mm]
\lesssim& m^{-2}\langle t\rangle^{-2}\sqrt{(\mathcal{E}_{12}+\mathcal{E}_{12}^{L})}\big(\langle t\rangle^{2}\mathcal{D}_{0}
+\langle t\rangle^{2}\mathcal{D}_{0}^{L}\big),\label{k31}\\[1mm]
J_4 \lesssim&
\|\eta\|_{3}\|\nabla q\|_0\|\eta^{d}\|_{2}\nonumber\\[1mm]
\lesssim&\|(\eta,\eta^{d})\|_{3}^2(\|u\|_{2}\|u\|_{3}+\|m\partial_{\omega}\eta\|_{4}\|m\partial_{\omega}\eta\|_{5})\nonumber\\[1mm]
\lesssim& m^{-2}\|m(\partial_{\omega}\eta,\partial_{\omega}\eta^{L})\|_{12}\|m(\partial_{\omega}\eta,\partial_{\omega}\eta^{L})\|_0
(\|u\|_{2}\|u\|_{3}
+\|m\partial_{\omega}\eta\|_{4}\|m\partial_{\omega}\eta\|_{5})\nonumber\\[1mm]
\lesssim& m^{-2}\langle t\rangle^{-2}{(\mathcal{E}_{12}+\mathcal{E}_{12}^{L})}\big(\langle t\rangle^{2}\mathcal{D}_{0}+\langle t\rangle^{2}\mathcal{D}_{0}^{L}\big).\label{k41}
\end{align}

\textbf{\emph{(2) Case $i\neq0$.}}

Integrating by part over $\Omega$ and using \eqref{poincareg}, \eqref{product21}, \eqref{al1}, \eqref{qestimate1} and the interpolation inequality, we can estimate that
\begin{align}
J_1 &=-\int\partial^{\alpha^{-}}\mm{curl}\mathcal{N}_{u}\cdot
\partial^{\alpha^{+}}\mm{curl}u^{d}\mm{d}y\lesssim\|\tilde{\ml{A}}(\tilde{\ml{A}}+I)\nabla u\|_{i+1}\|\mm{curl} u^{d}\|_{i+1} \nonumber\\[1mm]
&\lesssim(\|\eta\|_{3}\|u\|_{i+2}+\|\eta\|_{i+2}\|u\|_{3})\|\mm{curl}u^{d}\|_{i+1}\nonumber\\[2mm]
&\lesssim
\begin{cases}
m^{-1}(\|m\partial_{\omega}\eta\|_{3}^{1/2}\|u\|_{3}^{1/2}\|m\partial_{\omega}\eta\|_{9}^{1/2}\|u\|_{9}^{1/2}
+\|m\partial_{\omega}\eta\|_{9}\|u\|_{3})\|\mm{curl}u^{d}\|_{5}\\
\quad\lesssim C_{\epsilon}m^{-2}\langle t\rangle^{-1}(\langle t\rangle\mathcal{E}_{3})\mathcal{D}_{12}+\epsilon\|\mm{curl}u^{d}\|_5^2
& \;\hbox{ for }\; i=4;\\[2mm]
 m^{-1}\big(\|m\partial_{\omega}\eta\|_{6}\|u\|_9
+\|m\partial_{\omega}\eta\|_{12}\|u\|_3\big)\|\mm{curl}u^{d}\|_{8}\\
\quad\lesssim C_{\epsilon}m^{-2}\mathcal{E}_{12}\mathcal{D}_{12}+\epsilon\|\mm{curl}u^{d}\|_{8}^2
& \;\hbox{ for }\; i=7;\\[2mm]
 \big(m^{-1}\|m\partial_{\omega}\eta\|_{6}\|u\|_{10}
+m^{-3/4}\|m\partial_{\omega}\eta\|_{12}^{3/4}\|\eta\|_{13}^{1/4}\|u\|_3\big)\|\mm{curl}u^{d}\|_{9}\\
\quad\lesssim C_{\epsilon}\max\{m^{-3/2},m^{-2}\}\mathcal{E}_{12}\mathcal{D}_{12}
+\epsilon\|\mm{curl}u^{d}\|_{9}^2
& \;\hbox{ for }\; i=8.
\end{cases}\label{k12}
\end{align}
and
\begin{align}
J_2 &=
\int\partial^{\alpha^{-}}\mm{curl}(\nabla_{\tilde{\ml{A}}}q)\cdot\partial^{\alpha^{+}}\mm{curl}u^{d}\mm{d}y
\lesssim \|\tilde{\ml{A}}\nabla q\|_{i}\|\mm{curl}u^{d}\|_{i+1} \nonumber\\[1mm]
&\lesssim(\|\eta\|_{3}\| q\|_{i+1}+\|\eta\|_{i+1}\| q\|_{3})\|\mm{curl}u^{d}\|_{i+1}
\lesssim\|\eta\|_{i+1}\| q\|_{i+1}\|\mm{curl}u^{d}\|_{i+1}\nonumber\\[1mm]
&\lesssim\|\eta\|_{i+1}\big(\|u\|_{2}\|u\|_{i+1}+\|m\partial_{\omega}\eta\|_{6}\|m\partial_{\omega}\eta\|_{i+1}
+m\|\eta\|_{i}\|m\partial_{\omega}\eta\|_{4}
\big)\|\mm{curl}u^{d}\|_{i+1}\nonumber\\[1mm]
&\lesssim
\begin{cases}
 m^{-1}\|m\partial_{\omega}\eta\|_{8}\big(\|u\|_{2}\|u\|_{5}+\|m\partial_{\omega}\eta\|_{8}\|m\partial_{\omega}\eta\|_{3}
\big)\|\mm{curl}u^{d}\|_{5}\\
\quad\lesssim C_{\epsilon}m^{-2}\langle t\rangle^{-1}(\mathcal{E}_{12}(\langle t\rangle\mathcal{E}_{3}))\mathcal{D}_{12}+\epsilon\|\mm{curl}u^{d}\|_5^2
& \;\hbox{ for }\; i=4;\\[2mm]
m^{-1}\|m\partial_{\omega}\eta\|_{11}\big(\|u\|_{2}\|u\|_{8}+\|m\partial_{\omega}\eta\|_{10}\|m\partial_{\omega}\eta\|_{4}
\big)
\|\mm{curl}u^{d}\|_{8}\\
\quad\lesssim C_{\epsilon}m^{-2}\mathcal{E}_{12}\mathcal{E}_{4}\mathcal{D}_{12} +\epsilon\|\mm{curl}u^{d}\|_8^2
& \;\hbox{ for }\; i=7;\\[1mm]
m^{-1}\|m\partial_{\omega}\eta\|_{12}\big(\|u\|_{2}\|u\|_{9}+\|m\partial_{\omega}\eta\|_{11}\|m\partial_{\omega}\eta\|_{4}
\big)\|\mm{curl}u^{d}\|_{9}\\
\quad\lesssim C_{\epsilon} m^{-2} \mathcal{E}_{12}\mathcal{E}_{4}\mathcal{D}_{12} +\epsilon\|\mm{curl}u^{d}\|_9^2
& \;\hbox{ for }\; i=8.
\end{cases}\label{k22}
\end{align}

In the same manner,
\begin{align}
J_3 & = -\int\partial^{\alpha^{-}}\mm{curl}\mathcal{N}^{u}\cdot
\partial^{\alpha^{+}}\mm{curl}\eta^{d}\mm{d}y
\lesssim\|\tilde{\ml{A}}(\tilde{\ml{A}}+I)\nabla u\|_{i+1}\|\nabla \eta^{d}\|_{i+1} \nonumber\\[1mm]
&\lesssim(\|\eta\|_{3}\|u\|_{i+2}+\|\eta\|_{i+2}\|u\|_{3})\|\eta^{d}\|_{i+2}\nonumber\\[1mm]
&\lesssim
\begin{cases}
\|\eta\|_{6}\|\eta^{d}\|_{6}\|u\|_{5}
 & \\
\quad\lesssim
m^{-2}\|m\partial_{\omega}(\eta,\eta^{L})\|_{12}^{3/2}\|u\|\|_{10}^{1/2}\|u\|\|_{0}^{1/2}
\|m\partial_{\omega}(\eta,\eta^{L})\|_{0}^{1/2}
\\
\quad\lesssim m^{-2}\langle t\rangle^{-1}
\sqrt{\langle t\rangle^{2}(\mathcal{E}_{0}+\mathcal{E}_{0}^{L})}(\mathcal{D}_{12}+\mathcal{D}^{L}_{12})
& \;\hbox{ for }\; i=4;\\[1mm]
\|\eta\|_{9}\|\eta^{d}\|_{9}\|u\|_{8}\lesssim
m^{-2}\|m\partial_{\omega}(\eta,\eta^{L})\|_{12}^2\|u\|_{8}\\
\quad\lesssim m^{-2}\sqrt{\mathcal{E}_{12}}(\mathcal{D}_{12}+\mathcal{D}^{L}_{12})
& \;\hbox{ for }\; i=7;\\[1mm]
(\|\eta\|_{3}\|u\|_{10}+\|\eta\|_{10}\|u\|_{3})\|\eta^{d}\|_{10}\\
\quad\lesssim
m^{-7/4}\|m\partial_{\omega}\eta\|_{6}\|m(\partial_{\omega}\eta,\partial_{\omega}\eta^{L})\|_{12}^{3/4}
\|(\eta,\eta^{L})\|_{13}^{1/4}\|u\|_{10}\\
\quad\;\;+m^{-3/2}\|m(\partial_{\omega}\eta,\partial_{\omega}\eta^{L})\|_{12}^{3/2}\|(\eta,\eta^{L})\|_{13}^{1/2}\|u\|_3\\
\quad\lesssim \max\{m^{-3/2},m^{-7/4}\}
\sqrt{(\mathcal{E}_{12}+\mathcal{E}_{12}^{L})}(\mathcal{D}_{12}+\mathcal{D}^{L}_{12})
& \;\hbox{ for }\; i=8
\end{cases}\label{k32}
\end{align}
and
\begin{align}
J_4 &=
\int\partial^{\alpha^{-}}\mm{curl}(\nabla_{\tilde{\ml{A}}}q)\cdot\partial^{\alpha^{+}}\mm{curl}\eta^{d}\mm{d}y
\lesssim \|\tilde{\ml{A}}\nabla q\|_{i}\|\nabla \eta^{d}\|_{i+1} \nonumber\\[1mm]
&\lesssim(\|\eta\|_{3}\| q\|_{i+1}+\|\eta\|_{i+1}\| q\|_{3})\|\eta^{d}\|_{i+2}
\lesssim\|\eta\|_{i+1}\|\eta^{d}\|_{i+2}\| q\|_{i+1}\nonumber\\[1mm]
&\lesssim
\begin{cases}
m^{-2}\|m\partial_{\omega}\eta\|_{4}^{1/2}
\|m(\partial_{\omega}\eta,m\partial_{\omega}\eta^{L})\|_{12}^{3/2}
\|(u,m\partial_{\omega}\eta)\|_{4}^{3/2}\|(u,m\partial_{\omega}\eta)\|_{10}^{1/2}\\
\quad
\lesssim m^{-2}\langle t\rangle^{-1}({\langle t\rangle\mathcal{E}_{4}})(\mathcal{D}_{12}+\mathcal{D}^{L}_{12})
& \;\hbox{ for }\; i=4;\\[1mm]
 m^{-2}\|m(\partial_{\omega}\eta,\partial_{\omega}\eta^{L})\|_{12}^2
\big(\|u\|_{2}\|u\|_{8}+\|m\partial_{\omega}\eta\|_{4}\|m\partial_{\omega}\eta\|_{10}\big)\\
\quad\lesssim m^{-2}{\mathcal{E}_{12}}(\mathcal{D}_{12}+\mathcal{D}^{L}_{12})
& \;\hbox{ for }\; i=7;\\[1mm]
m^{-7/4}\|m\partial_{\omega}\eta\|_{12}\|m(\partial_{\omega}\eta,\partial_{\omega}\eta^{L})\|_{12}^{3/4}
\|(\eta,\eta^{L})\|_{13}^{1/4}\|q\|_{9}\\
\quad\lesssim m^{-7/4}
(\mathcal{E}_{12}+\mathcal{E}_{12}^{L})\mathcal{D}_{12}
& \;\hbox{ for }\; i=8.
\end{cases}\label{k42}
\end{align}
{
Therefore, it follows from \eqref{202210151940} and \eqref{k11}--\eqref{k42} that,
for a given suitable small $\epsilon$,
\begin{align}
&\frac{\mm{d}}{\mm{d}t}\mathcal{E}^{d}_{\mm{curl},i}+\mathcal{D}^{d}_{\mm{curl},i}\nonumber\\
&\lesssim
\begin{cases}
m^{-2}\big(\sqrt{\mathcal{E}_{12}}+\mathcal{E}_{12}+\mathcal{E}_{12}^2)(\mathcal{D}_{12}+\mathcal{D}^{L}_{12}\big)
& \;\hbox{ for }\; i=7;\\[1mm]
\max\{m^{-3/2}, m^{-7/4},m^{-2}\}
\big(\sqrt{(\mathcal{E}_{12}+\mathcal{E}_{12}^{L})}+\mathcal{E}_{12}+\mathcal{E}_{12}^{L}+\mathcal{E}_{12}^2\big)
(\mathcal{D}_{12}+\mathcal{D}^{L}_{12})
& \;\hbox{ for }\; i=8
\end{cases}\label{202210151621}
\end{align}
and
\begin{align}\label{202210152012}
&\sum_{j=0}^{2}\langle m^{-1}\rangle^{2j}\left(\frac{\mm{d}}{\mm{d}t}\langle t\rangle^{(2-j)}\mathcal{E}_{\mm{curl},4j}^{d}
+\langle t\rangle^{(2-j)}\mathcal{D}^{d}_{\mm{curl},4j}\right)\nonumber\\
&\lesssim m^{-2}\bigg((\langle t\rangle^{2}\mathcal{E}_{0})(1+\mathcal{E}_{12})\mathcal{D}_{12}
+\big(\sqrt{(\mathcal{E}_{12}+\mathcal{E}_{12}^{L})}+{\mathcal{E}_{12}+\mathcal{E}_{12}^{L}}\big)\big(\langle t\rangle^{2}\mathcal{D}_{0}+\langle t\rangle^{2}\mathcal{D}_{0}^{L}\big)\bigg)\nonumber\\
&\;\;\;+\langle m^{-1}\rangle^2 m^{-2}(1+\mathcal{E}_{12})
\bigg(\sqrt{\langle t\rangle^{2}(\mathcal{E}_{0}+\mathcal{E}_{0}^{L})}
+{\langle t\rangle\mathcal{E}_{4}}\bigg)\big(\mathcal{D}_{12}+\mathcal{D}^{L}_{12}\big)\nonumber\\
&\;\;\;+\langle m^{-1}\rangle^4\max\{m^{-3/2}, m^{-7/4},m^{-2}\}
\bigg(\sqrt{(\mathcal{E}_{12}+\mathcal{E}_{12}^{L})}+\mathcal{E}_{12}+\mathcal{E}_{12}^{L}+\mathcal{E}_{12}^2\bigg)
\big(\mathcal{D}_{12}+\mathcal{D}^{L}_{12}\big).
\end{align}
Integrating  \eqref{202210151621} and \eqref{202210152012} with respect to time over $(0,t)$, and then utilizing \eqref{202207041516n}--\eqref{202207041516nn}, \eqref{202109090810}, and \eqref{20221014112n}--\eqref{202210141123}, we further deduce that
\begin{align}
&\mathcal{E}^{d}_{\mm{curl},i}(t)+\int_{0}^{t}\mathcal{D}^{d}_{\mm{curl},i}(\tau)\mm{d}\tau\nonumber\\[1mm]
&\lesssim\|(\nabla\mm{curl}\eta^{d}, \mm{curl}u^{d}, m\mm{curl}\partial_{\omega}\eta^{d})(0)\|_{i}^2\nonumber\\[1mm]
&\quad+
\begin{cases}
m^{-2}\big(\sqrt{\ml{E}_{\mm{total}}^0}+\ml{E}_{\mm{total}}^0\big)^3
& \;\hbox{ for }\; i=7;\\[2mm]
\max\{m^{-3/2}, m^{-7/4},m^{-2}\}\big(\sqrt{\ml{E}_{\mm{total}}^0}+\ml{E}_{\mm{total}}^0\big)^3
& \;\hbox{ for }\; i=8
\end{cases}\label{202210151621nn}
\end{align}
and
\begin{align}\label{202210152012nn}
&\sum_{j=0}^{2}\langle m^{-1}\rangle^{2j}\left(\langle t\rangle^{(2-j)}\mathcal{E}_{\mm{curl},4j}^{d}(t)
+\int_{0}^{t}\langle \tau\rangle^{(2-j)}\mathcal{D}^{d}_{\mm{curl},4j}(\tau)\mm{d}\tau\right)\nonumber\\
&\lesssim \sum_{j=0}^{2}\langle m^{-1}\rangle^{2j}\|(\nabla\mm{curl}\eta^{d}, \mm{curl}u^{d}, m\mm{curl}\partial_{\omega}\eta^{d})(0)\|_{4j}^2\nonumber\\
&\lesssim m^{-2}\big(\sqrt{\ml{E}_{\mm{total}}^0}+\ml{E}_{\mm{total}}^0\big)^2\Xi
+\langle m^{-1}\rangle^2m^{-2}\big(\sqrt{\ml{E}_{\mm{total}}^0}+\ml{E}_{\mm{total}}^0\big)^2\Xi\nonumber\\[1mm]
&\;\;\;+\langle m^{-1}\rangle^4\max\{m^{-3/2}, m^{-7/4},m^{-2}\}
\big(\sqrt{\ml{E}_{\mm{total}}^0}+\ml{E}_{\mm{total}}^0\big)^3\nonumber\\[1mm]
&\lesssim \langle m^{-1}\rangle^2m^{-2}\big(\sqrt{\ml{E}_{\mm{total}}^0}+\ml{E}_{\mm{total}}^0\big)^2\Xi\nonumber\\[1mm]
&\quad
+\langle m^{-1}\rangle^4\max\{m^{-3/2}, m^{-7/4},m^{-2}\}\big(\sqrt{\ml{E}_{\mm{total}}^0}+\ml{E}_{\mm{total}}^0\big)^3.
\end{align}}

Additionally,
recalling that $(\eta^{d})_{\mathbb{T}^3}=(u^{d})_{\mathbb{T}^3}=0$, $\mm{div}u^{d}=-\mm{div}_{\tilde{\ml{A}}}u$ and $\mm{div}\eta^{d}=\mm{div}\eta$ for any $t>0$, we can have
$$(\nabla\mm{div}\eta^{d}, \mm{div}u^{d}, m\mm{div}\partial_{\omega}\eta^{d})
=(\nabla\mm{div}\eta, -\mm{div}_{\tilde{\ml{A}}}u, m\mm{div}\partial_{\omega}\eta).$$
By making use of \eqref{202209210921}--\eqref{202209210921nn}, \eqref{poincareg} and the interpolation inequality, we find that
\begin{align}\label{202210161400}
&\|(\nabla\mm{div}\eta^{d}, \mm{div}u^{d}, m\mm{div}\partial_{\omega}\eta^{d})\|_{i}^2\nonumber\\
&\lesssim
\begin{cases}
\|\eta\|_{3}^2\|(\nabla\eta, u, m\partial_{\omega}\eta)\|_{1}^2& \;\hbox{ for }\; i=0;\\[1mm]
\|\eta\|_{3}^2\|(\nabla\eta, u, m\partial_{\omega}\eta)\|_{i+1}^2
+\|\eta\|_{i+1}^2\|(\nabla\eta, u, m\partial_{\omega}\eta)\|_{3}^2& \;\hbox{ for }\; i=4,\ 7,\ 8
\end{cases}\nonumber\\
&\lesssim
\begin{cases}
m^{-2}\|m\partial_{\omega}\eta\|_{6}^2\|(\nabla\eta, u, m\partial_{\omega}\eta)\|_{1}^2& \;\hbox{ for }\; i=0;\\[1mm]
m^{-2}\|m\partial_{\omega}\eta\|_{6}^2\|(\nabla\eta, u, m\partial_{\omega}\eta)\|_{i+1}^2
+m^{-2}\|m\partial_{\omega}\eta\|_{i+4}^2\|(\nabla\eta, u, m\partial_{\omega}\eta)\|_{3}^2& \;\hbox{ for }\; i=4,\ 7,\ 8
\end{cases}\nonumber\\
&\lesssim
\begin{cases}
m^{-2}\langle t\rangle^{-2}(\langle t\rangle^{2}\mathcal{E}_{0})\mathcal{E}_{12}& \;\hbox{ for }\; i=0;\\[1mm]
m^{-2}\langle t\rangle^{-1}(\langle t\rangle\mathcal{E}_{4})\mathcal{E}_{8}& \;\hbox{ for }\; i=4;\\[1mm]
m^{-2}\mathcal{E}_{12}\mathcal{E}_{6}  & \;\hbox{ for }\; i=7,\ 8
\end{cases}
\end{align}
and
\begin{align}\label{202210161402}
&\|(\nabla\mm{div}u^{d}, m\mm{div}\partial_{\omega}\eta^{d})\|_{i}^2\nonumber\\
&\lesssim
\begin{cases}
\|\eta\|_{3}^2\|(\nabla u, m\partial_{\omega}\eta)\|_{i+1}^2& \;\hbox{ for }\; i=0;\\[1mm]
\|\eta\|_{4}^2\|(\nabla u, m\partial_{\omega}\eta)\|_{i+1}^2
+\|\eta\|_{i+2}^2\|(\nabla u, m\partial_{\omega}\eta)\|_{3}^2& \;\hbox{ for }\; i=4,7,\ 8
\end{cases}\nonumber\\
&\lesssim
\begin{cases}
m^{-2}\|m\partial_{\omega}\eta\|_{6}^2\|(\nabla u, m\partial_{\omega}\eta)\|_{i+1}^2& \;\hbox{ for }\; i=0;\\[1mm]
m^{-2}\|m\partial_{\omega}\eta\|_{7}^2\|(\nabla u, m\partial_{\omega}\eta)\|_{i+1}^2
+m^{-2}\|m\partial_{\omega}\eta\|_{i+5}^2\|(u, m\partial_{\omega}\eta)\|_{3}^2& \;\hbox{ for }\; i=4,\ 7
;\\[1mm]
m^{-2}\|m\partial_{\omega}\eta\|_{7}^2\|(\nabla u, m\partial_{\omega}\eta)\|_{i+1}^2
+m^{-3/2}\|m\partial_{\omega}\eta\|_{12}^{3/2}\|\eta\|_{13}^{1/2}
\|(u, m\partial_{\omega}\eta)\|_{3}^2& \;\hbox{ for }i=8
\end{cases}\nonumber\\
&\lesssim
\begin{cases}
m^{-2}\langle t\rangle^{-2}(\langle t\rangle^{2}\mathcal{D}_{0})\mathcal{E}_{12}& \;\hbox{ for }\; i=0;\\[1mm]
m^{-2}\langle t\rangle^{-1}(\langle t\rangle\mathcal{E}_{4})\mathcal{D}_{12}& \;\hbox{ for }\; i=4;\\[1mm]
m^{-2}\mathcal{E}_{12}\mathcal{D}_{4}& \;\hbox{ for }\; i=7;\\[1mm]
\max\{m^{-2},m^{-3/2}\}\mathcal{E}_{12}\mathcal{D}_{6}& \;\hbox{ for }\; i=8.
\end{cases}
\end{align}
Combining \eqref{202210161400} with \eqref{202210161402} then gives rise to
\begin{align}
&\|(\nabla\mm{div}\eta^{d}, \mm{div}u^{d}, m\mm{div}\partial_{\omega}\eta^{d})(t)\|_{i}^2
+\int_{0}^{t}\|(\nabla\mm{div}u^{d}, m\mm{div}\partial_{\omega}\eta^{d})(\tau)\|_{i}^2\mm{d}\tau\nonumber\\[1mm]
&\lesssim
\begin{cases}
m^{-2}\big(\ml{E}_{\mm{total}}^0\big)^2
& \;\hbox{ for }\; i=7;\\[2mm]
\max\{m^{-2},m^{-3/2}\}\ml{E}_{\mm{total}}^0\ml{E}_{6}^0
& \;\hbox{ for }\; i=8
\end{cases}\label{202210151621nnnn}
\end{align}
and
\begin{align}\label{202210152012nnnnn}
&\sum_{j=0}^{2}\bigg(\langle m^{-1}\rangle^{2j}\langle t\rangle^{(2-j)}\|(\nabla\mm{div}\eta^{d}, \mm{div}u^{d}, m\mm{div}\partial_{\omega}\eta^{d})(t)\|_{4j}^2\nonumber\\
&\quad\quad\quad+\langle m^{-1}\rangle^{2j}\int_{0}^{t}\langle \tau\rangle^{(2-j)}\|(\nabla\mm{div}u^{d}, m\mm{div}\partial_{\omega}\eta^{d})(\tau)\|_{4j}^2\mm{d}\tau\bigg)\nonumber\\[1mm]
&\lesssim
m^{-2}\ml{E}_{\mm{total}}^0\Xi+\langle m^{-1}\rangle^2m^{-2}\ml{E}_{\mm{total}}^0\Xi
+\langle m^{-1}\rangle^4
\max\{m^{-3/2},m^{-2}\}\ml{E}_{\mm{total}}^0\Xi\nonumber\\[2mm]
&\lesssim
\langle m^{-1}\rangle^2m^{-2}\ml{E}_{\mm{total}}^0\Xi
+\langle m^{-1}\rangle^4
\max\{m^{-3/2},m^{-2}\}\ml{E}_{\mm{total}}^0\ml{E}_{6}^0.
\end{align}

Furthermore, by virtue of \eqref{202210161627} and the interpolation inequality, we have that
\begin{align}
\label{202210291950}
&\|(\nabla\mm{curl}\eta^{d}, \mm{curl}u^{d}, m\mm{curl}\partial_{\omega}\eta^{d})(0)\|_{i}^2\nonumber
\\
&\lesssim
\|(\nabla\eta^{\mm{r}}, u^{\mm{r}}, m\partial_{\omega}\eta^{\mm{r}})\|_{i+1}^2\nonumber\\&
\lesssim
\begin{cases}
m^{-2}\mathcal{E}_{6}^0\mathcal{E}_{2}^0
& \;\hbox{ for }\; i=0;\\[2mm]
m^{-2}\big(\mathcal{E}_{6}^0\mathcal{E}_{5}^0+\mathcal{E}_{8}^0\mathcal{E}_{3}^0\big)\lesssim m^{-2}\mathcal{E}_{8}^0\mathcal{E}_{3}^0
& \;\hbox{ for }\; i=4;\\[2mm]
m^{-2}\big(\mathcal{E}_{6}^0\mathcal{E}_{8}^0+\mathcal{E}_{12}^0\mathcal{E}_{3}^0\big)\lesssim m^{-2}\mathcal{E}_{12}^0\mathcal{E}_{3}^0
& \;\hbox{ for }\; i=7;\\[1mm]
m^{-2}\big(\mathcal{E}_{6}^0\mathcal{E}_{9}^0+\mathcal{E}_{12}^0\mathcal{E}_{3}^0\big)\lesssim m^{-2}\mathcal{E}_{12}^0\mathcal{E}_{3}^0
& \;\hbox{ for }\; i=8.
\end{cases}
\end{align}
Consequently,
with \eqref{202210161627},  \eqref{202210151621nn}--\eqref{202210152012}, \eqref{202210151621nnnn}--\eqref{202210152012nnnnn} and \eqref{202210291950} in hand, and recalling
$(\eta^{d})_{\mathbb{T}^3}=(u^{d})_{\mathbb{T}^3}=0$, we then utilize \eqref{hodgeellipticnn} to obtain
\begin{align*}
&\mathcal{E}^{d}_{i+1}(t)+\int_{0}^{t}\mathcal{D}^{d}_{i+1}(\tau)\mm{d}\tau\nonumber\\
&\lesssim
\begin{cases}
m^{-2}\big(\sqrt{\ml{E}_{\mm{total}}^0}+\ml{E}_{\mm{total}}^0\big)^3
& \;\hbox{ for }\; i=7;\\[2mm]
\max\{m^{-3/2}, m^{-7/4},m^{-2}\}\big(\sqrt{\ml{E}_{\mm{total}}^0}+\ml{E}_{\mm{total}}^0\big)^3
& \;\hbox{ for }\; i=8
\end{cases}
\end{align*}
and
\begin{align*}
&\sum_{j=0}^{2}\langle m^{-1}\rangle^{2j}\left(\langle t\rangle^{(2-j)}\mathcal{E}_{4j+1}^{d}(t)
+\int_{0}^{t}\langle \tau\rangle^{(2-j)}\mathcal{D}^{d}_{4j+1}(\tau)\mm{d}\tau\right)\nonumber\\
&\lesssim
\langle m^{-1}\rangle^2m^{-2}\big(\sqrt{\ml{E}_{\mm{total}}^0}+\ml{E}_{\mm{total}}^0\big)^2\Xi
+\langle m^{-1}\rangle^4\max\{m^{-3/2}, m^{-7/4},m^{-2}\}\big(\sqrt{\ml{E}_{\mm{total}}^0}+\ml{E}_{\mm{total}}^0\big)^3.
\end{align*}
The above two inequalities
then yields \eqref{202210011616nnnn}--\eqref{202210011616nn} immediately.
This completes the proof of Theorem \ref{thm:2022nn}.

\vspace{4mm}
\noindent\textbf{Acknowledgements.}
The research of Youyi Zhao was supported by Fujian Alliance Of Mathematics (2023SXLMQN01) and NSFC (Grant No.  12371233). The author is deeply grateful to Prof. Fei Jiang for kind guidance. In addition, the author thanks Boyang Fang for carefully checking this paper.

%
%

\renewcommand\refname{References}
\renewenvironment{thebibliography}[1]{%
\section*{\refname}
\list{{\arabic{enumi}}}{\def\makelabel##1{\hss{##1}}\topsep=0mm
\parsep=0mm
\partopsep=0mm\itemsep=0mm
\labelsep=1ex\itemindent=0mm
\settowidth\labelwidth{\small[#1]}%
\leftmargin\labelwidth \advance\leftmargin\labelsep
\advance\leftmargin -\itemindent
\usecounter{enumi}}\small
\def\newblock{\ }
\sloppy\clubpenalty4000\widowpenalty4000
\sfcode`\.=1000\relax}{\endlist}
\bibliographystyle{model1b-num-names}

\begin{thebibliography}{22}
\expandafter\ifx\csname natexlab\endcsname\relax\def\natexlab#1{#1}\fi
\providecommand{\bibinfo}[2]{#2}
\ifx\xfnm\relax \def\xfnm[#1]{\unskip,\space#1}\fi
\bibitem[{Adams and John(2005)}]{ARAJJFF}
\bibinfo{author}{R.A. Adams}, \bibinfo{author}{J.F.F. John},
  \bibinfo{title}{{Sobolev Space}}, \bibinfo{publisher}{Academic Press: New
  York}, \bibinfo{year}{2005}.
\bibitem[{Bardos et~al.(1988)Bardos, Sulem and Sulem}]{BCSCSPLL}
\bibinfo{author}{C.~Bardos}, \bibinfo{author}{C.~Sulem},
  \bibinfo{author}{P.~Sulem}, \bibinfo{title}{{Longtime dynamics of a
  conductive fluid in the presence of a strong magnetic field}},
  \bibinfo{journal}{Trans. Amer. Math. Soc.} \bibinfo{volume}{305}
  (\bibinfo{year}{1988}) \bibinfo{pages}{175--191}.
\bibitem[{Califano and Chiuderi(1999)}]{CFCCRI}
\bibinfo{author}{F.~Califano}, \bibinfo{author}{C.~Chiuderi},
  \bibinfo{title}{Resistivity-independent dissipation of magnetohydrodynamic
  waves in an inhomogeneous plasma}, \bibinfo{journal}{Phy. Rev. E}
  \bibinfo{volume}{60} (\bibinfo{year}{1999}) \bibinfo{pages}{4701--4707}.
\bibitem[{Chandrasekhar(1961)}]{CSHHSCPO}
\bibinfo{author}{S.~Chandrasekhar}, \bibinfo{title}{{Hydrodynamic and
  Hydromagnetic Stability, The International Series of Monographs on Physics}},
  \bibinfo{publisher}{Oxford, Clarendon Press}, \bibinfo{year}{1961}.
\bibitem[{Chen et~al.(2021)Chen, Zhang and Zhou}]{CZZ2021}
\bibinfo{author}{W.J. Chen}, \bibinfo{author}{Z.F. Zhang},
  \bibinfo{author}{J.F. Zhou}, \bibinfo{title}{{Global well-posedness for the
  3-D MHD equations with partial diffusion in the periodic domain}},
  \bibinfo{journal}{Sci. China Math.} \bibinfo{volume}{64}
  (\bibinfo{year}{2021}) \bibinfo{pages}{1--10}.
\bibitem[{Fujita and Kato(1985)}]{FHKT}
\bibinfo{author}{H.~Fujita}, \bibinfo{author}{T.~Kato}, \bibinfo{title}{{On the
  Navier-Stokes initial value problem. I}}, \bibinfo{journal}{Arch. Ration.
  Mech. Anal.} \bibinfo{volume}{62} (\bibinfo{year}{1985})
  \bibinfo{pages}{167--186}.
\bibitem[{Guo and Tice(2013{\natexlab{a}})}]{GYTIAE2}
\bibinfo{author}{Y.~Guo}, \bibinfo{author}{I.~Tice}, \bibinfo{title}{{Almost
  exponential decay of periodic viscous surface waves without surface
  tension}}, \bibinfo{journal}{Arch. Ration. Mech. Anal.} \bibinfo{volume}{207}
  (\bibinfo{year}{2013}{\natexlab{a}}) \bibinfo{pages}{459--531}.
\bibitem[{Guo and Tice(2013{\natexlab{b}})}]{GYTIDAP}
\bibinfo{author}{Y.~Guo}, \bibinfo{author}{I.~Tice}, \bibinfo{title}{{Decay of
  viscous surface waves without surface tension in horizontally infinite
  domains}}, \bibinfo{journal}{Anal. PDE} \bibinfo{volume}{6}
  (\bibinfo{year}{2013}{\natexlab{b}}) \bibinfo{pages}{1429--1533}.
\bibitem[{Jiang and Jiang(2018)}]{JFJSJMFMOSERT}
\bibinfo{author}{F.~Jiang}, \bibinfo{author}{S.~Jiang}, \bibinfo{title}{{ On
  the stabilizing effect of the magnetic fields in the magnetic
  Rayleigh--Taylor problem}}, \bibinfo{journal}{SIAM J. Math. Anal.}
  \bibinfo{volume}{50} (\bibinfo{year}{2018}) \bibinfo{pages}{491--540}.
\bibitem[{Jiang and Jiang(2019{\natexlab{a}})}]{JFJSSETEFP}
\bibinfo{author}{F.~Jiang}, \bibinfo{author}{S.~Jiang}, \bibinfo{title}{{ On
  the dynamical stability and instability of Parker problem}},
  \bibinfo{journal}{Physica D} \bibinfo{volume}{391}
  (\bibinfo{year}{2019}{\natexlab{a}}) \bibinfo{pages}{17--51}.
\bibitem[{Jiang and Jiang(2019{\natexlab{b}})}]{JFJSCVPDE1}
\bibinfo{author}{F.~Jiang}, \bibinfo{author}{S.~Jiang},
  \bibinfo{title}{{Nonlinear stability and instability in the Rayleigh--Taylor
  problem of stratified compressible MHD fluids}}, \bibinfo{journal}{Calc. Var.
  Partial Differential Equations} \bibinfo{volume}{58}
  (\bibinfo{year}{2019}{\natexlab{b}}) \bibinfo{pages}{29}.
\bibitem[{Jiang and Jiang(2019{\natexlab{c}})}]{JFJSOMITIN}
\bibinfo{author}{F.~Jiang}, \bibinfo{author}{S.~Jiang}, \bibinfo{title}{{On
  magnetic inhibition theory in non-resistive magnetohydrodynamic fluids}},
  \bibinfo{journal}{Arch. Ration. Mech. Anal.} \bibinfo{volume}{233}
  (\bibinfo{year}{2019}{\natexlab{c}}) \bibinfo{pages}{749--798}.
\bibitem[{Jiang and Jiang(2020)}]{JFJSOUI}
\bibinfo{author}{F.~Jiang}, \bibinfo{author}{S.~Jiang}, \bibinfo{title}{{ On
  inhibition of thermal convection by a magnetic field under zero
  resistivity}}, \bibinfo{journal}{J. Math. Pures Appl.} \bibinfo{volume}{141}
  (\bibinfo{year}{2020}) \bibinfo{pages}{220--265}.
\bibitem[{Jiang and Jiang(2021)}]{JFJSO2020}
\bibinfo{author}{F.~Jiang}, \bibinfo{author}{S.~Jiang},
  \bibinfo{title}{{Asymptotic behaviors of global solutions to the
  two-dimensional non-resistive MHD equations with large initial
  perturbations}}, \bibinfo{journal}{Adv. Math.} \bibinfo{volume}{393}
  (\bibinfo{year}{2021}) \bibinfo{pages}{108084}.
\bibitem[{Jiang and Jiang(2023)}]{JFJSO2021}
\bibinfo{author}{F.~Jiang}, \bibinfo{author}{S.~Jiang}, \bibinfo{title}{{On
  magnetic inhibition theory in non-resistive magnetohydrodynamic fluids:
  global existence of large solutions}}, \bibinfo{journal}{Arch. Ration. Mech.
  Anal.} \bibinfo{volume}{247} (\bibinfo{year}{2023}) \bibinfo{pages}{96}.
\bibitem[{Kraichnan(1965)}]{RHK}
\bibinfo{author}{R.H. Kraichnan}, \bibinfo{title}{{Inertial-range spectrum of
  hrdromagnetic turbulence}}, \bibinfo{journal}{Phys. Fluids}
  \bibinfo{volume}{8} (\bibinfo{year}{1965}) \bibinfo{pages}{1385--1387}.
\bibitem[{Majda and Bertozzi(2002)}]{MAJBAL}
\bibinfo{author}{A.J. Majda}, \bibinfo{author}{A.L. Bertozzi},
  \bibinfo{title}{{Vorticity and incompressible flow}},
  \bibinfo{publisher}{Cambridge university press}, \bibinfo{year}{2002}.
\bibitem[{Novotn{\`y} and Stra{\v{s}}kraba(2004)}]{NASII04}
\bibinfo{author}{A.~Novotn{\`y}}, \bibinfo{author}{I.~Stra{\v{s}}kraba},
  \bibinfo{title}{{Introduction to the Mathematical Theory of Compressible
  Flow}}, \bibinfo{publisher}{Oxford University Press, USA},
  \bibinfo{year}{2004}.
\bibitem[{Ren et~al.(2014)Ren, Wu, Xiang and Zhang}]{RXXWJHXZYZZF}
\bibinfo{author}{X.X. Ren}, \bibinfo{author}{J.H. Wu}, \bibinfo{author}{Z.Y.
  Xiang}, \bibinfo{author}{Z.F. Zhang}, \bibinfo{title}{{Global existence and
  decay of smooth solution for the 2D MHD equations without magnetic
  diffusion}}, \bibinfo{journal}{J. Funct. Anal.} \bibinfo{volume}{267}
  (\bibinfo{year}{2014}) \bibinfo{pages}{503--541}.
\bibitem[{Tan and Wang(2018)}]{TZWYJGw}
\bibinfo{author}{Z.~Tan}, \bibinfo{author}{Y.J. Wang}, \bibinfo{title}{{Global
  well-posedness of an initial-boundary value problem for viscous non-resistive
  MHD systems}}, \bibinfo{journal}{SIAM J. Math. Anal.} \bibinfo{volume}{50}
  (\bibinfo{year}{2018}) \bibinfo{pages}{1432--1470}.
\bibitem[{Wang(2019)}]{WYTIVNMI}
\bibinfo{author}{Y.J. Wang}, \bibinfo{title}{{ Sharp nonlinear stability
  criterion of viscous non-resistive MHD internal waves in 3D}},
  \bibinfo{journal}{Arch. Ration. Mech. Anal.} \bibinfo{volume}{231}
  (\bibinfo{year}{2019}) \bibinfo{pages}{1675--1743}.
\bibitem[{Zhang(2016)}]{ZTGS}
\bibinfo{author}{T.~Zhang}, \bibinfo{title}{{Global solutions to the 2D
  viscous, non-resistive MHD system with large background magnetic field}},
  \bibinfo{journal}{J. Differential Equations} \bibinfo{volume}{260}
  (\bibinfo{year}{2016}) \bibinfo{pages}{5450--5480}.

\end{thebibliography}

\end{document}